\documentclass[a4paper,12pt]{amsart}
\usepackage{amsfonts}
\usepackage{amssymb}
\usepackage{ifthen}
\usepackage{graphicx}
\nonstopmode \numberwithin{equation}{section}
\usepackage[usenames]{color}
\newtheorem{thm}{Theorem}[section]
\newtheorem{lem}{Lemma}[section]
\newtheorem{cor}[thm]{Corollary}
\newtheorem{prop}{Proposition}[section]

\theoremstyle{definition}
\newtheorem{mlem}{Main lemma}[section]
\newtheorem{assertion}{Assertion}[section]
\newtheorem{cl}{Claim}[section]
\newtheorem{ca}{Case}[section]
\newtheorem{sca}{Subcase}[section]
\newtheorem{scl}{Subclaim}[section]
\newtheorem{conj}[thm]{Conjecture}
\newtheorem{fact}{Fact}[section]
\newtheorem{defn}{Definition}[section]
\newtheorem{op}[thm]{Open Problem}
\newtheorem{ques}[thm]{Question}
\newtheorem{rem}{Remark}[section]
\newtheorem{exam}{Example}[section]

\newtheorem{condition}[thm]{Condition}

\numberwithin{equation}{section}

\newcounter {own}
\def\theown {\thesection       .\arabic{own}}

\newenvironment{pf}[1][]{%
 \vskip 3mm
 \noindent
 \ifthenelse{\equal{#1}{}}%
  {{\slshape Proof. }}%
  {{\slshape #1.} }%
 }%
{\qed\bigskip}

\newcounter{alphabet}

\newenvironment{Thm}[1][]{\refstepcounter{alphabet}%
\bigskip%
\noindent%
{\bf Theorem \Alph{alphabet}}%
\ifthenelse{\equal{#1}{}}{}{ (#1)}%
{\bf .} \itshape}{\vskip 8pt}

\newcounter{alphabet2}

\newcounter{alphabet3}

\newenvironment{Cons}[1][]{\refstepcounter{alphabet3}%
\bigskip%
\noindent%
{\bf Condition \Alph{alphabet3}}%
\ifthenelse{\equal{#1}{}}{}{ (#1)}%
{\bf .} \itshape}{\vskip 8pt}

\makeatletter
\renewcommand{\Ref}[1]{\ref{#1}}
\makeatother

\newenvironment{Lem}[1][]{\refstepcounter{alphabet}%
\bigskip%
\noindent%
{\bf Lemma \Alph{alphabet}}%
{\bf .} \itshape}{\vskip 8pt}

\newcommand{\IR}{{\mathbb R}}

\newcommand{\IB}{{\mathbb B}}

\newcommand{\diam}{{\operatorname{diam}}}

\newcommand{\dist}{{\operatorname{dist}}}

\def\be{\begin{equation}}
\def\ee{\end{equation}}

\newcommand{\ben}{\begin{enumerate}}
\newcommand{\een}{\end{enumerate}}

\newcommand{\blem}{\begin{lem}}
\newcommand{\elem}{\end{lem}}
\newcommand{\bthm}{\begin{thm}}
\newcommand{\ethm}{\end{thm}}
\newcommand{\bcor}{\begin{cor}}
\newcommand{\ecor}{\end{cor}}
\newcommand{\beg}{\begin{exam}}
\newcommand{\eeg}{\end{exam}}
\newcommand{\begs}{\begin{examples}}
\newcommand{\eegs}{\end{examples}}
\newcommand{\bdefe}{\begin{defn}}
\newcommand{\edefe}{\end{defn}}
\newcommand{\bprob}{\begin{prob}}
\newcommand{\eprob}{\end{prob}}
\newcommand{\bques}{\begin{ques}}
\newcommand{\eques}{\end{ques}}
\newcommand{\bei}{\begin{itemize}}
\newcommand{\eei}{\end{itemize}}
\newcommand{\bcon}{\begin{conj}}
\newcommand{\econ}{\end{conj}}
\newcommand{\bop}{\begin{op}}
\newcommand{\eop}{\end{op}}

\newcommand{\bas}{\begin{assertion}}
\newcommand{\eas}{\end{assertion}}

\newcommand{\bfa}{\begin{fact}}
\newcommand{\efa}{\end{fact}}

\newcommand{\bca}{\begin{ca}}
\newcommand{\eca}{\end{ca}}
\newcommand{\bsca}{\begin{sca}}
\newcommand{\esca}{\end{sca}}

\newcommand{\bcl}{\begin{cl}}
\newcommand{\ecl}{\end{cl}}

\newcommand{\bmlem}{\begin{mlem}}
\newcommand{\emlem}{\end{mlem}}

\newcommand{\bscl}{\begin{scl}}
\newcommand{\escl}{\end{scl}}

\newcommand{\bcons}{\begin{conjs}}
\newcommand{\econs}{\end{conjs}}

\newcommand{\bprop}{\begin{prop}}
\newcommand{\eprop}{\end{prop}}

\newcommand{\bcond}{\begin{condition}}
\newcommand{\econd}{\end{condition}}

\newcommand{\br}{\begin{rem}}
\newcommand{\er}{\end{rem}}
\newcommand{\brs}{\begin{rems}}
\newcommand{\ers}{\end{rems}}
\newcommand{\bo}{\begin{obser}}
\newcommand{\eo}{\end{obser}}
\newcommand{\bos}{\begin{obsers}}
\newcommand{\eos}{\end{obsers}}
\newcommand{\bpf}{\begin{pf}}
\newcommand{\epf}{\end{pf}}
\newcommand{\ba}{\begin{array}}
\newcommand{\ea}{\end{array}}
\newcommand{\beq}{\begin{eqnarray}}
\newcommand{\beqq}{\begin{eqnarray*}}
\newcommand{\eeq}{\end{eqnarray}}
\newcommand{\eeqq}{\end{eqnarray*}}

\newcounter{minutes}
\divide\time by 60
\newcounter{hours}
\multiply\time by 60 \addtocounter{minutes}{-\time}
\begin{document}

\bibliographystyle{amsplain}

\title{Geometric characterizations of inner uniformity through Gromov hyperbolicity}

\author{Manzi Huang}
\address{Manzi Huang, MOE-LCSM, School of Mathematics and Statistics, Hunan Normal University, Changsha, Hunan 410081, People's Republic of China} \email{mzhuang@hunnu.edu.cn}

\author{Antti Rasila}
\address{Antti Rasila,
Mathematics with Computer Science Program, Guangdong Technion, 241 Daxue Road, Jinping District, Shantou, Guangdong 515063, People's Republic of China and Department of Mathematics, Technion - Israel Institute of Technology, Haifa 32000, Israel}
\email{antti.rasila@gtiit.edu.cn; antti.rasila@technion.ac.il; antti.rasila@iki.fi}

\author{Xiantao Wang}
\address{Xiantao Wang, MOE-LCSM, School of Mathematics and Statistics, Hunan Normal University, Changsha, Hunan 410081, People's Republic of China}
\email{xtwang@hunnu.edu.cn}

\author{Qingshan Zhou}%${}^{~\mathbf{*}}$}
\address{Qingshan Zhou,
%Department of Mathematics, Shantou University, Shantou, Guangdong 515063, People's Republic of China}
  School of Mathematics, Foshan University,  Foshan, Guangdong 528000, People's Republic
of China}
\email{q476308142@qq.com}

\date{}
\keywords{Geometric characterization, inner uniformity, Gromov hyperbolicity, Gromov boundary, quasisymmetry, linear local connectedness.
%\\ ${}^{\mathbf{*}}$ Corresponding author
}

\begin{abstract}
In this paper, we study the characterization of inner uniformity of bounded domains $G$ in $\IR^n$, and prove that the following three conditions are equivalent: $(1)$ $G$ is inner uniform;
$(2)$ $G$ is Gromov hyperbolic and its inner metric boundary is naturally quasisymmetrically equivalent to the Gromov boundary; $(3)$ $G$ is Gromov hyperbolic and linearly locally connected with respect to the inner metric.
The equivalence between the conditions $(1)$ and $(2)$, and the implication from $(2)$ to $(3)$ affirmatively answer three questions raised by Bonk, Heinonen, and Koskela in 2001.
\end{abstract}

\thanks{}

\maketitle \pagestyle{myheadings} \markboth{Manzi Huang, Antti Rasila, Xiantao Wang and Qingshan Zhou}{Geometric characterizations of inner uniformity through Gromov hyperbolicity}

\maketitle

%%%%%%%    \tableofcontents

%%%%%%%%%%%%%%%%%%%%%%%%%%%%%%%

%%%%%%%%%%%%%%%%%%%%%%%%%%%%%%%%%%%%%%%%%%%
%%%%%%%%%%%%%%%%%%%%%%%%%%%%%%%%%%%%%%%%%%%
\section{Introduction} \label{sec-1}
%%%%%%%%%%%%%%%%%%%%%%%%%%%%%%%%%%%%%%%%%%%
%%%%%%%%%%%%%%%%%%%%%%%%%%%%%%%%%%%%%%%%%%%

Gromov hyperbolicity is a concept introduced by Gromov in the setting of geometric group theory in 1980s \cite{Gr87}.
Since its introduction, Gromov hyperbolicity has found numerous applications, and it has been, for example,
considered in the books \cite{BH99, BBI01, CDP90, GH90, Ro03, Sh91}. Initially, it was mainly studied from the point of view of the hyperbolic group theory (see, e.g., \cite{GH90}). Recently,
geometric characterizations of Gromov hyperbolicity have been established in \cite{BB03, KLM14}, the connection between Gromov hyperbolicity and quasiconformal deformations has been studied in \cite{Her06}, and the Gromov hyperbolicity of various metrics and surfaces has been investigated in  \cite{Bo06, Ha06, MM99, RT06}. For other discussions in this line, see, for example, \cite{ BS, BK15, BH20,HL21, HLV07, Vai10, We08}.

A domain, i.e., an open and connected set, in a metric space is said to be {\it uniform} if every pair of points can be joined by a so called uniform curve (see Section \ref{sec-2} for the precise definition).
Uniform domains were independently introduced by John in \cite{Jo61} and by Martio and Sarvas in \cite{MS78}. Since its introduction, this concept has played a significant role in the study of geometric function theory, see \cite{HRWZ} and references therein. Also, uniform domains have been generalized to Carnot groups, including Heisenberg groups \cite{CT95, CGN00, Gre01}, as well as the general metric spaces \cite{BSh-07, FOS-15, LPZ-22, OSh-17, Raj-21}.

In \cite{BHK}, Bonk, Heinonen, and Koskela considered the relationship between Gromov hyperbolicity and uniformity of domains in $\IR^n$, and established the following characterization.

\begin{Thm}\label{Thm-BHK}  $($\cite[Theorem 1.11]{BHK}$)$
A bounded domain in $\IR^n$ is uniform if and only if $(a)$ it is Gromov hyperbolic, and $(b)$ its Euclidean boundary is naturally quasisymmetrically equivalent to the Gromov boundary.
\end{Thm}
See also \cite{HSX, Vai05} for its analogues in Banach spaces and metric spaces, respectively.

The terminology in Theorem \Ref{Thm-BHK} and in the rest of this section will be introduced in the second section unless otherwise stated.

If the metric of a metric space is replaced by the corresponding inner metric, then uniformity is changed into inner uniformity.
Inner uniform domains were studied in the plane by Balogh and Volberg in their study of the complex iteration of certain polynomials \cite{BV96-1, BV96-2}, where they called these domains uniformly John domains. For a comprehensive survey on inner uniformity and related concepts, see \cite{Vai058}.
Recently, there has been substantial interest in study of inner uniformity. For example, the invariance of inner uniformity under quasiconformal mappings in $\IR^n$ has been investigated in \cite{BuSt2}. A boundary Harnack principle for inner uniform domains has been established
in \cite{Lie15, LS14-1, LS14-2}. Neumann and Dirichlet heat kernels on inner uniform domains have been considered in \cite{GS11, Sal10}. The characterization of Gromov hyperbolic domains on the $2$-sphere in terms of inner uniformity, and the relationships between inner uniformity and certain capacity conditions have been studied in \cite{BHK, BH07, BH}. In particular, in \cite{BHK}, Bonk, Heinonen, and Koskela proved that inner uniformity implies Gromov hyperbolicity (see \cite[Theorem 1.11]{BHK}).
A natural question is whether there is a result similar to the characterization of uniform domains in Theorem \Ref{Thm-BHK} for inner uniform domains.
In fact, Bonk, Heinonen, and Koskela raised this as an open problem (see the second paragraph after Theorem $1.11$ in \cite{BHK}). For convenience, we call this open problem the {\it Characterization Question} in the following.

We remark that, in the proof of Theorem $7.11$ in \cite{BHK}, i.e., the sufficiency part of Theorem \Ref{Thm-BHK}, the authors first proved that the linear local connectedness of Gromov hyperbolic domains implies their uniformity \cite[Proposition 7.12]{BHK}, and then, verified that the condition $(b)$ in Theorem \Ref{Thm-BHK} ensures the linear local connectedness of such domains \cite[Proposition 7.13]{BHK}. In the same paper, the authors
proposed the questions whether \cite[Theorem 7.11]{BHK} and \cite[Proposition 7.13]{BHK} are true when the Euclidean metric is replaced by the corresponding inner metric (see the second paragraph after Theorem $7.11$ and the second paragraph after Proposition $7.13$, respectively, in \cite{BHK}). For convenience, we call the first open problem the {\it Sufficiency Question} and the second one the {\it LLC Question}. Obviously, an affirmative answer to  the Characterization Question will imply that the answer to the Sufficiency Question is also affirmative.

Recent studies, in particular, ones related to the abovementioned three questions, indicate that the relationships of Gromov hyperbolicity, inner uniformity and linear local connectedness of domains in $\IR^n$ are crucial topics  of investigation. The purpose of this paper is to study these topics. Our result is as follows.

\begin{thm}\label{thm-1.1} Suppose that $G$ is a bounded domain in $\IR^n$ with $n\geq 2$. Then the following conditions are equivalent:
\begin{enumerate}
\item[$(\mathfrak{i})$]\label{thm-1.1-1}
$G$ is inner uniform;

\item[$(\mathfrak{ii})$]\label{thm-1.1-2}
$G$ is Gromov hyperbolic and its inner metric boundary is naturally quasisymmetrically equivalent to the Gromov boundary;

\item[$(\mathfrak{iii})$]\label{thm-1.1-3}
$G$ is Gromov hyperbolic and {\it LLC} $($the abbreviation of ``linearly locally connected"$)$ with respect to the inner metric $\sigma$.
\end{enumerate}
\end{thm}

Let $(X,d)$ denote a metric space, and let $c\geq 1$ be a constant. A domain $G\subsetneq X$ is called {\it $c$-linearly locally connected}, or $c$-{\it LLC}, with respect to the metric $d$ if the following conditions hold (cf. \cite{BHK} or \cite{BH}): For all $x\in G$ and all $r>0$,
\begin{enumerate}
\item[$LLC_1$:]
every pair of points in $\mathbb{B}_d(x,r)\cap G$ can be joined by a rectifiable curve in $\mathbb{B}_d(x,c r)\cap G$.
\item[$LLC_2$:]
every pair of points in $G\setminus \overline{\mathbb{B}}_d(x,r)$ can be joined by a rectifiable curve in  $G\setminus \overline{\mathbb{B}}_d(x,r/c)$.
\end{enumerate}
Here, $\mathbb{B}_d(x,r)$ denotes the metric ball in $X$ with center $x$ and radius $r$, i.e., $\mathbb{B}_d(x,r):=\{y\in X:\; d(y,x)<r\}$, and $\overline{\mathbb{B}}_d(x,r):=\{y\in X:\; d(y,x)\leq r\}$.

The concept of linear local connectedness was first introduced by Gehring when he discussed the characterizations of quasidisks \cite{FW87}, and it is well-known in the literature.

As a direct consequence of the equivalence of the conditions $(\mathfrak{ii})$ and $(\mathfrak{iii})$ in Theorem \ref{thm-1.1}, we see that for a bounded Gromov hyperbolic domain in $\IR^n$, the natural quasisymmetrical equivalence of its inner metric boundary and the Gromov boundary is equivalent to the {\it LLC} property with respect to the inner metric. Consequently, it follows that the answer to the {\it LLC} Question is affirmative.

In the proof of Theorem \ref{thm-1.1}, we use an approach that is completely different from the ideas applied in the proof of \cite[Theorem 1.11]{BHK}. The outline of the proof is as follows. First, we check the equivalence of the conditions $(\mathfrak{i})$ and $(\mathfrak{ii})$ in Theorem \ref{thm-1.1}.
The implication from  $(\mathfrak{i})$ to  $(\mathfrak{ii})$, denoted by $(\mathfrak{i})$ $\Longrightarrow$ $(\mathfrak{ii})$, easily follows from \cite[Theorem 1.31]{BHK} and \cite[Theorem $2.39$]{Vai05} (i.e., Theorem \Ref{sat-12} below). When we verify the opposite implication, i.e., the one $(\mathfrak{ii})$ $\Longrightarrow$ $(\mathfrak{i})$, difficulties arise from the fact that no conditions related to mappings acting in $G$ are assumed in $(\mathfrak{ii})$, as the inner uniformity is a property defined in terms of the interior geometry of $G$.
To overcome this difficulty, we introduce Condition \ref{PropertyB} below, and then, prove the implication $(\mathfrak{ii})$ $\Longrightarrow$ $(\mathfrak{i})$
by showing the ones:  $(\mathfrak{ii})$ $\Longrightarrow$ Condition \ref{PropertyB} $\Longrightarrow$  $(\mathfrak{i})$. Thus the needed equivalence is proved. This check is carried out in Sections \ref{sec-2}$-$\ref{sec-6}.
 This equivalence also implies that the answer to the Characterization Question is affirmative, and hence, the answer to
the Sufficiency Question is affirmative as well.

Based on the equivalence of $(\mathfrak{i})$ and $(\mathfrak{ii})$, combining with Balogh and Buckley's geometric characterization of Gromov hyperbolicity, i.e., Theorem \Ref{BB03} below, in Section \ref{sec-7}, we prove the implications:  $(\mathfrak{ii})$ $\Longrightarrow$  $(\mathfrak{iii})$ $\Longrightarrow$  $(\mathfrak{i})$, and thus, the proof of Theorem \ref{thm-1.1} is complete.

In Section \ref{sec-2}, necessary definitions and terminology will be presented, useful known results will be recalled, and two lemmas will be proved.

%%%%%%%%%%%%%%%%%%%%%%%%%%%%%%%%%%%%%%%%%%%
%%%%%%%%%%%%%%%%%%%%%%%%%%%%%%%%%%%%%%%%%%%
\section{Preliminaries}\label{sec-2}
%%%%%%%%%%%%%%%%%%%%%%%%%%%%%%%%%%%%%%%%%%%
%%%%%%%%%%%%%%%%%%%%%%%%%%%%%%%%%%%%%%%%%%%

Let $(X, d)$ denote a metric space. A {\it curve} in $X$ is a continuous function
$\gamma:\; I\to X$ from an interval $I\subset \IR$ to $X$. If $\gamma$ is an embedding of $I$,
it is also called an {\it arc}. We use $\gamma$ to denote both the function and its image set. The {\it length} $\ell_d(\gamma)$ of $\gamma$ with respect to the metric $d$ is defined in the usual way. The parameter interval $I$ is allowed to be closed, open or half-open. If $\ell_d(\gamma) < \infty$, then $\gamma$ is said to be {\it rectifiable}.

The metric space $(X, d)$ is called {\it rectifiably connected} if every pair of points in $X$ can be joined by a rectifiable
curve in $X$, and {\it geodesic} if every pair of points $x,$ $y$ in $X$ can be joined
by a curve $\gamma$ with $\ell_d(\gamma)=d(x, y)$.

For convenience, in the following, we assume that $(X, d)$, $(Y, d')$ denote rectifiably connected, locally compact and non-complete metric spaces.

\subsection{Uniform domains, inner uniform domains and John domains}
Let $M\geq 1$ be a constant, and let $\gamma:$ $[0, 1]\to X$ denote a curve with endpoints $x=\gamma(0)$ and $y=\gamma(1)$. For points $z,w$ on $\gamma$, let $\gamma[z,w]$ be the restriction of $\gamma$ to the interval $[t_0,t_1]$, where $t_0=\min\{t_z,t_w\}$, $t_1=\max\{t_z,t_w\}$, $t_z =\max\{t\in [0,1]: \gamma(t)=z\}$, and $t_w =\min\{t\in [0,1]: \gamma(t)=w\}$.

Let $D\subset X$ be a domain. For a curve $\gamma\subset D$ with endpoints $x$ and $y$ in $D$, we say that $\gamma$ satisfies:
\ben
\item
 the {\it $M$-John property} (with respect to $d$-length) if
$$\min\{\ell_d(\gamma[x,z]),\ell_d(\gamma[z,y])\}\leq M\delta(d)_{D,X}(z)$$ holds for any $z=\gamma(t)$, where $t\in [0,1]$ and $\delta(d)_{D,X}(z)$ denotes the distance from $z$ to the boundary of $D$ with respect to the metric $d$.
\item
the {\it $M$-quasiconvexity } if $\ell_d(\gamma)\leq Md(x, y)$.
%\item
%{\it $M$-turning condition $($diameter$)$} if $\diam(\gamma)\leq Md(x, y)$.
\een

Let $\overline{X_d}$ denote the metric completion $($with respect to the metric $d$$)$ and
$$\partial_dX=\overline{X_d}\setminus X,$$ which is called the {\it $($metric$)$ boundary} of $X$.
For a rectifiably connected, locally compact and non-complete metric space $(X, d)$, we know that $\partial_dX$ is nonempty and closed, that is, $\delta(d)_X(x)>0$ for any $x\in X$ (cf. \cite[2.A]{BH}).
For such spaces, the definitions for John property and quasiconvexity extend, in an obvious way, to the situations where the
parameter interval is open or half open.

A domain $D$ in $X$ is called \ben
\item
{\it $M$-John} if any two points in $D$ can be joined by a curve in $D$ satisfying the $M$-John property;
\item
{\it $M$-uniform} if any two points in $D$ can be joined by an $M$-uniform curve in $D$, where a curve is called {\it $M$-uniform} if it satisfies both the $M$-John property and the $M$-quasiconvexity.
\een

The {\it inner metric} $\sigma(d)$ of $d$ is defined as follows: For any pair of points $x$ and $y$ in $D$,
$$\sigma(d)(x, y)=\inf\{\ell_d(\gamma)\},$$
where the infimum is taken over all rectifiable curves $\gamma$ in $D$ connecting $x$ and $y$.

If $D$ is $M$-uniform with respect to $\sigma(d)$, then $D$ is called {\it inner $M$-uniform}. Similarly, we may define the concept of {\it inner $M$-uniform curves}.

\br  Every $M$-uniform domain is inner $M$-uniform, but inner uniformity does not imply uniformity (cf. \cite[Page3]{BHK}). The following relations follow immediately from the definitions: An inner $M$-uniform domain is $M$-John, $M$-John implies $M_1$-John when $M\leq M_1$. Similarly, (inner) $M$-uniformity of domains implies (inner) $M_1$-uniformity provided that $M\leq M_1$.
\er
\blem\label{changsha-2}
For a domain $D$ in $X$, suppose that $D$ is $M$-John with $M\geq 1$. Then $D$ is bounded with respect to the metric $d$, if and only if it is bounded with respect to the inner metric $\sigma(d)$.
\elem
\bpf
First we come to prove the sufficiency. For any $x$, $y\in D$, we know from the definition of $\sigma(d)$ that
$$d(x,y)\leq \sigma(d)(x,y),$$
and so $$\diam_{d}(D)\leq \diam_{\sigma(d)}(D).$$ This implies that the sufficiency is true.

Next, we will prove the necessity. For this, we let $x$, $y\in D$. Then there is a curve $\gamma$ in $D$ connecting $x$ and $y$ such that for any $z\in \gamma$,
$$\min\{\ell_d(\gamma[x,z]),\ell_d(\gamma[z,y])\}\leq M\delta(d)_{D,X}(z).$$
Let $z_0\in \gamma$ bisect $\gamma$, i.e., $\ell_d(\gamma[x,z_0])=\ell_d(\gamma[y,z_0])$. Then we have
$$\ell_d(\gamma)=2\ell_d(\gamma[x,z_0])\leq 2M\delta(d)_{D,X}(z_0)\leq 2M\diam_d(D),$$ from which the necessity part follows,
where $\diam_d(D)$ denotes the diameter of $D$ with respect to $d$.
\epf

\br  We know from \cite[Lemma 1.3]{HRWZ} that there exists a domain $D$ which is bounded with respect to the Euclidean metric $|\cdot|$, whereas, it is not bounded with respect to the inner metric $\sigma(|\cdot|)$. This shows that the assumption that the domain $D$ is John in Lemma \ref{changsha-2} cannot be removed.
\er

\subsection{Quasihyperbolic metrics}
The {\it quasihyperbolic length} of a rectifiable curve
$\gamma$ in a rectifiably connected, locally compact and non-complete metric space $(X, d)$ is the number:
$$\ell_{k(d)_X}(\gamma)=\int_{\gamma}\frac{|dz|}{\delta(d)_X(z)},
$$ where $|dz|$ denotes the length element in $X$ with respect to the metric $d$.

For any $x$, $y$ in $X$,
the {\it quasihyperbolic distance}
$k(d)_X(x,y)$  between $x$ and $y$ is defined by
$$k(d)_X(x,y)=\inf\{\ell_{k(d)_X}(\gamma)\},
$$ where the infimum is taken over all rectifiable curves $\gamma$ in $X$ with endpoints $x$ and $y$. The resulting metric space $(X, k(d)_X)$ is complete, proper and geodesic provided that the identity mapping $\mathop{\mathrm{id}}:$ $(X,d)\to (X,\sigma(d))$ is a homeomorphism (cf. \cite{BHK}).

For a rectifiable curve $\gamma$ in $X$ connecting $x$ and $y$, the following estimate on $\ell_{k(d)_X}(\gamma)$ is useful (cf. \cite{BHK}):
\beq\label{vai-lem3.2}\ell_{k(d)_X}(\gamma)\geq \log\Big(1+\frac{\ell_d(\gamma)}
{\min\{\delta(d)_X(x), \delta(d)_X(y)\}}\Big),
\eeq
and thus,
\beq\label{eq-neq-1}
k(d)_X(x,y)\geq \log\Big(1+\frac{\sigma(d)(x,y)}
{\min\{\delta(d)_X(x), \delta(d)_X(y)\}}\Big)\geq \Big|\log \frac{\delta(d)_X(x)}{\delta(d)_X(y)}\Big|.
\eeq

Recall that a curve $\gamma$ connecting $x$ and
$y$ is a {\it quasihyperbolic geodesic} if $\ell_{k(d)_X}(\gamma)=k(d)_X(x,y)$.  Each subcurve of a quasihyperbolic
geodesic is a quasihyperbolic geodesic. It is known that every proper subdomain $D$ of $\IR^n$ $(n\geq 2)$ is quasihyperbolic geodesic connected, see \cite[Lemma 1]{GO79}. This is not true in arbitrary metric spaces (see \cite[Theorem 2.3]{RT-14} or \cite[Example 2.9]{Vai90}).

For other basic properties of the quasihyperbolic metrics, the reader is referred to \cite{GO79, GP76}.

%%%%%%%%%%%%%%%%%%%%%%%%%%%%%%%%%%%%%%%%%%%

\subsection{Quasigeodesics and quasigeodesic rays}
For a given constant $\lambda\geq 1$, a rectifiable curve or a rectifiable ray $\gamma$ in $X$ is called {\it $\lambda$-quasigeodesic} if for any two points $u$ and $v$ in $\gamma$,
$$\ell_{k(d)_X}(\gamma[u,v])\leq \lambda k(d)_X(u,v),$$ where {\it a ray} in $X$ means a curve with one of its endpoints in $X$ and the other on $\partial_d X$.

Obviously, a $\lambda$-quasigeodesic (resp. a $\lambda$-quasigeodesic ray) is a quasihyperbolic geodesic (resp. a quasihyperbolically geodesic ray) if and only if $\lambda=1$.

In 1991, V\"ais\"al\"a established the following property concerning the existence
of quasigeodesics in Banach spaces: Suppose that $D$ is a proper subdomain in a Banach space and $\lambda> 1$ is a constant. Then for any pair of points in $D$, there is a
$\lambda$-quasigeodesic in $D$ joining these two points (see \cite[Theorem 3.3]{Vai91}). Here and in the following, all Banach spaces are assumed to have dimension at least $2$.

\subsection{Conformal deformations, Gromov hyperbolic domains and Gromov hyperbolic spaces}
Let us recall the following conformal deformations (cf. \cite[Chapter $4$]{BHK}).  Fix a base point $p\in X$, and consider the family of conformal deformations of $X$ defined by the densities
$$\rho(d)_{\varepsilon,p}(x)=e^{-\varepsilon d(x,p)}\;\;(\varepsilon>0).$$

For $u$, $v\in X$, let
$$(d)_{\varepsilon,p}(u,v)=\inf\int_{\gamma} \rho(d)_{\varepsilon,p}(x) ds,$$
where the infimum is taken over all rectifiable curves $\gamma$ in $X$ connecting $u$ and $v$.
Then $(d)_{\varepsilon,p}$ are metrics on $X$.
% We denote the resulting metric spaces by $X_\varepsilon=(X,(d)_\varepsilon)$.

 Suppose that $\delta\geq 0$ is a constant. \ben
 \item\label{def-1(1)}
 We say that $(X, d)$ is {\it Gromov $\delta$-hyperbolic} if
for all $x,$ $y,$ $z,$ $p\in X$,
$$(x|y)_p \geq \min\{(x|z)_p, (z|y)_p\}-\delta,$$
 where $(x|y)_p$ is the {\it Gromov product} defined by
$$ 2(x|y)_p =d(x,p)+d(y,p)-d(x,y).$$
\item\label{def-1(2)}
For a proper subdomain $D$ in $X$, it is called {\it Gromov $\delta$-hyperbolic} if
$(D,k(d)_D)$ is Gromov $\delta$-hyperbolic.
\een

Also, we say that a metric space is {\it Gromov hyperbolic} if it is Gromov $\delta$-hyperbolic for some $\delta\geq 0$. It is known that all (inner) uniform domains in $\IR^n$ are Gromov hyperbolic (see \cite[Theorem 1.11]{BHK}).

We remark that the Gromov hyperbolicity as defined in \eqref{def-1(1)} is equivalent to the one given below in geodesic metric spaces (cf. \cite{BS07}).

Let $(X, d)$ be geodesic and $\delta$ a nonnegative constant. Denote by
$[x, y]$ any geodesic joining two points $x$ and $y$ in $X$. For $x_1$, $x_2$, $x_3\in D$, let $\alpha_1=[x_1,x_2]$, $\alpha_2=[x_2,x_3]$ and $\alpha_3=[x_1,x_3]$, and we denote $\Delta=(\alpha_1, \alpha_2, \alpha_3)$ to be the geodesic triangle. We say the  geodesic triangle $\Delta$ is {\it $\delta$-thin } if 
$$d(w, \alpha_{i+1}\cup \alpha_{i+2})\leq \delta$$ holds for any $w\in \alpha_i$ with $i\in\{1,2,3\}$, where 
$\alpha_4=\alpha_1$ and $\alpha_5=\alpha_2$.
 $(X, d)$ is called {\it Gromov $\delta$-hyperbolic} if every geodesic triangle in $X$ is {\it $\delta$-thin.}

The following theorem says that the deformations $X_\varepsilon$ are uniform whenever $(X,d)$ is a proper, geodesic and Gromov hyperbolic space.

\begin{Thm}\label{Lem-1}$($\cite[Proposition $4.5$]{BHK}$)$
There is a constant $\varepsilon_0=\varepsilon_0(\delta)>0$ such that for any $0<\varepsilon<\varepsilon_0$ and $p\in X$, the conformal deformation $(X,(d)_{\varepsilon,p})$ of a proper, geodesic and Gromov $\delta$-hyperbolic space $(X,d)$ is a bounded $A(\delta)$-uniform space, where the notation $\varepsilon_0(\delta)$ $($resp. $A(\delta)$$)$ means that the constant
$\varepsilon_0$ $($resp. $A$$)$ depends only on $\delta$.
\end{Thm}

Suppose that $(X, d)$ is Gromov $\delta$-hyperbolic and $p\in X$ is a base point.\ben
\item
A sequence $\{x_i\}$ in $X$ is called a {\it Gromov sequence} if $(x_i|x_j)_p\to \infty$ as $i,$ $j\to \infty.$
\item
Two Gromov sequences $\{x_i\}$ and $\{y_j\}$ are said to be {\it equivalent} if $(x_i|y_i)_p\to \infty$ as $i\to \infty.$
\item
The {\it Gromov boundary} $\partial^*X$ of $X$ is defined to be the set of all equivalent classes, and $X^*=X \cup \partial^*X$ is called the {\it Gromov closure} of $X$.
\item
For $a\in X$ and $b\in \partial^*X$, the Gromov product $(a|b)_p$ of $a$ and $b$ is defined by
$$(a|b)_p= \inf \big\{ \liminf_{i\to \infty}(a|b_i)_p:\; \{b_i\}\in b\big\}.$$
\item
For $a,$ $b\in \partial^*X$, the Gromov product $(a|b)_p$ of $a$ and $b$ is defined by
$$(a|b)_p= \inf \big\{ \liminf_{i\to \infty}(a_i|b_i)_p:\; \{a_i\}\in a\;\;{\rm and}\;\; \{b_i\}\in b\big\}.$$
\een

\subsection{Geometric characterization of Gromov hyperbolicity}
We say that a domain $G\subsetneq \mathbb{R}^n$ satisfies (cf. \cite{BB03}, \cite{BHK} or \cite{KLM14})
\begin{enumerate}
\item the {\it Gehring-Hayman condition} if there is a constant $c_{gh}\geq 1$ such that for any pair of points $x$, $y\in G$, for any quasihyperbolic geodesic in $G$ connecting $x$ and $y$, denoted by $[x, y]_k$, and for every curve $\gamma$ in $G$ with endpoints $x$ and $y$, the following inequality holds:
    $$\ell([x, y]_k)\leq c_{gh}\ell(\gamma),$$
    where $\ell([x, y]_k)$ (resp. $\ell(\gamma)$) stands for the Euclidean length of $[x, y]_k$ (resp. $\gamma$);
\item the {\it ball-separation condition} if there is a constant $c_{bs}\geq 1$ such that for any pair of points $x$, $y\in G$, for any quasihyperbolic geodesic $[x, y]_k$, for any point $z\in [x, y]_k$, and for every curve $\gamma$ in $G$ with endpoints $x$ and $y$, the following relation holds:
    $$\mathbb{B}_\sigma(z, c_{bs}\delta_G(z))\cap \gamma \neq \emptyset,
    $$ where $\delta_G(z)$ means the distance from $z$ to the boundary of $G$ (with respect to the Euclidean metric).
\end{enumerate}

The following characterization of Gromov hyperbolicity of proper subdomains in $\IR^n$ is useful for our discussions (cf. \cite[Theorem $0.1$]{BB03} or \cite[Theorem 1.1]{KLM14}).

\begin{Thm}\label{BB03}
Suppose that $G\varsubsetneq \IR^n$ is a domain. Then $G$ is Gromov $\delta$-hyperbolic if and only if it satisfies both the $c_{gh}$-Gehring-Hayman condition and the $c_{bs}$-ball-separation condition. The constants $\delta\geq 0$, $c_{gh}\geq 1$ and $c_{bs}\geq 1$ depend on each other, and on $n$.
\end{Thm}

\subsection{Visual metrics and natural mappings}\label{changsha-1}
Suppose that $(X, d)$ is a Gromov $\delta$-hyperbolic space. For $p\in X$ and $\tau>0$, define $$\rho_{p,\tau}(x,y)=e^{-\tau(x|y)_p}$$ for $x,$ $y\in X^*$ with convention $e^{-\infty}=0$.
Then it follows from \cite[Proposition $5.16$]{Vai10} (see also \cite[\S3]{BHK}) that there is a constant $\tau_0=\tau_0(\delta)>0$ such that for any $0<\tau<\tau_0$, one can define
a function $(d)_{\tau}^p$ which satisfies
\beqq\label{eq-11}(d)_{\tau}^p\leq \rho_{p,\tau}\leq 2(d)_{\tau}^p,\eeqq
where the function $(d)_{\tau}^p$ is a metametric on $X^*$, that is, it satisfies the axioms of a metric except that $(d)_{\tau}^p(x,x)$ may be positive. In fact, $(d)_{\tau}^p(x,y)=0$ if and only if $x=y\in \partial^*X$. Hence $(d)_{\tau}^p$ defines a metric on $\partial^*X$, which is called the {\it visual metric} of $\partial^*X$.

The metametric $(d)_{\tau}^p$ defines a topology $\mathcal{T}^*$ in $X^*$. In this topology, the points of $X$ are isolated. For a sequence $\{x_i\}\in X$ and $a\in \partial^*X$, $(d)_{\tau}^p(a, x_i)\to 0$ as $i\to\infty$ if and only if $\{x_i\}$ is a Gromov sequence and $\{x_i\}\in a$ (see \cite[Lemma 5.3]{Vai10}).

For a domain $D$ in $X$, let $D^{\star}=(D,k(d)_D)\cup \partial^*(D,k(d)_D)$. Since the restriction $\mathcal{T}^*|_{D}$ is discrete, the identity mapping $\mathop{\mathrm{id}}:$ $D \to D$ is continuous from the topology $\mathcal{T}^*$ to the metric topology of $D$. If it has a continuous extension
$$\varphi:\;\;(D^*, (k(d)_D)^p_\tau)\to(\overline{D_d},d),$$ then we call $\varphi$ a {\it natural mapping}.

Suppose that $E$ is a Banach space with metric $d$ and dimension at least $2$. The following two results, due to V\"{a}is\"{a}l\"{a}, are useful:

\begin{Thm} $($\cite[Lemma $2.22$]{Vai05}$)$\label{qz-1}
Suppose that $D\subset E$ is a Gromov $\delta$-hyperbolic domain. Then the natural mapping $$\varphi:\;\;(D^*, (k(d)_D)^p_\tau)\to(\overline{D_d},d)$$ exists if and only if every Gromov sequence $\overline{x}=\{x_i\}$ in $(D, k(d)_D)$ has a limit $\xi$ with respect to $d$. Moreover,  for each $\eta\in \partial^* D$ and for any Gromov sequence $\overline{x}\in \eta$, $\overline{x}$ converges to $\xi\in \partial_d D$ in $(\overline{D_d},d)$ and $\xi=\varphi(\eta)$.
\end{Thm}

\begin{Thm} $($\cite[Proposition $2.26$]{Vai05}$)$\label{qz-2}
Suppose that $D\subset E$ is an $M$-uniform domain. Then the natural mapping
$$\varphi:\;\;(D^*, (k(d)_D)^p_\tau)\to(\overline{D_d},d)$$ exists and is bijective for any $\tau$ with $0<\tau\leq \min\{1, \tau_0\}$. Moreover, a sequence $\overline{x}=\{x_i\}$ in $D$ converges to $\xi\in \partial_d D$ with respect to $d$ if and only if $\overline{x}$ is a Gromov sequence in $(D, k(d)_D)$ and $\varphi(\eta)=\xi\in \partial_d D$, where $\eta\in \partial^* D$ with $\overline{x}\in \eta$.
\end{Thm}

%\medskip

\subsection{Useful classes of mappings}\label{sub-2.5}

Let $f:$ $(X,d)\to (Y,d')$ be a mapping (not necessarily continuous), and let $M\geq 1$ and $K\geq 0$ be constants. If $\sup_{w\in Y}\{d'(w, f(X))\}<+\infty$ and for all $x,y\in X$,
$$M^{-1}d(x,y)-K\leq d'(f(x),f(y))\leq Md(x,y)+K,$$
then $f$ is called an {\it $(M, K)$-roughly quasi-isometric mapping} (cf. \cite{BS} or \cite{HSX}), where $d'(w, f(X))$ denotes the distance from the point $w$ to the image $f(X)$ of $X$ under $f$ with respect to $d'$.
An $(M, 0)$-roughly quasi-isometric mapping is said to be {\it $M$-quasi-isometric}.
If we replace $(X,d)$ by $(I,|\cdot|)$, where $I$ denotes an interval in $\IR$, then $f$ is called an {\it $(M, K)$-roughly quasi-isometric curve} (cf. \cite{Vai10}).

Suppose that both $(X, d)$ and $(Y, d')$ are rectifiably connected, locally compact and non-complete.
A homeomorphism $f:$ $(X,d)\to (Y,d')$ is called {\it $M$-quasihyperbolic} with $M\geq 1$ if
$$M^{-1}k(d)_X(x,y)\leq k(d')_{Y}(f(x),f(y))\leq M k(d)_X(x,y)$$
for all $x$, $y\in X$.

Obviously, a homeomorphism between two proper subdomains in metric spaces is $M$-quasihyperbolic if and only if it is $M$-quasi-isometric (or bi-Lipschitz) with respect to the corresponding quasihyperbolic metrics.

Suppose that $\eta$ is a self-homeomorphism of $[0, \infty)$.  A homeomorphism $f:$ $(X,d)\to (Y,d')$ is said to be
{\it $\eta$-quasisymmetric} if $d(a,x)\leq t d(x,b)$ implies
\beq\label{quasisymmetry} d'(f(a),f(x))\leq \eta(t) d'(f(x),f(b))\eeq for all $t\geq 0$ and for each triple $\{a,x,b\}$ in $X$.

Quasisymmetric mappings originate from the work of Beurling and Ahlfors \cite{BA56}, who defined them as the boundary values of quasiconformal self-mappings of the upper half-plane on the real line. The definition of quasisymmetric mappings above is due to
Tukia and V\"ais\"al\"a, who introduced the general class of quasisymmetric mappings \cite{TV80}. Since its introduction, the concept has been generalized to metametrics (see \cite[Subsection 4.3]{Vai10}), and has been studied by numerous authors, see \cite{HRWZ} and the references therein.

If we assume that the inequality \eqref{quasisymmetry} in the definition of quasisymmetric mappings holds for each triple $\{a,x,b\}$ in $X$ with $x\in A$ or $\{a, b\}\subset A$, then $f$ is called {\it $\eta$-quasisymmetric rel $A$}.

The following two theorems due to V\"{a}is\"{a}l\"{a} are used in the discussions in Section \ref{sec-3}.

\begin{Thm}\label{sat-12} $($\cite[Theorem $2.39$]{Vai05}$)$
Suppose that $D$ is a bounded $M$-uniform domain in a Banach space $E$ with a
base point $p\in D$ such that
$$\delta(d)_D(x)\leq c\delta(d)_D(p)$$
for all $x\in D$, where $d$ denotes the norm metric in $E$ and $c\geq 1$ is a constant. Then the bijective natural
mapping $$\psi:\;\;(D^*,(k(d)_D)^p_{\mu})\to (\overline{D_d},d),$$ which exists by Theorem \Ref{qz-2}, is $\eta$-quasisymmetric
rel $\partial^*D$ with respect to the metametric $(k(d)_D)^p_{\mu}$ of $D^*$ and the metric $d$ of $\overline{D_d}$, where $0< \mu\leq \mu_0,$ $\eta=\eta_{c, M, \mu}$ $($i.e., the control function $\eta$ depends only on the given parameters $c,$ $M$ and $\mu)$, and $\mu_0=\mu_0(M)$.
\end{Thm}

\begin{Thm}\label{fri-2} $($\cite[Theorem $5.35$]{Vai10}$)$
Suppose that $(X,p)$ and $(Y,p')$ are pointed length Gromov $\delta$-hyperbolic spaces and that $f:$ $X\to Y$ is a $(\lambda,\mu)$-roughly quasi-isometric mapping with $p'=f(p)$.
Then $f$ has an extension
$$f^*:\;\;(X^*,(d)^p_\varepsilon)\to (Y^*,(d')^{p'}_\varepsilon),$$ which is continuous, where $0<\varepsilon \leq \min\{1,\tau_0\}$, $d$ and $d'$ denote the metrics in $X$ and $Y$, respectively. Moreover,
 if $f$ is weakly surjective, then
the restriction $$f^*|_{\partial^*X}:\;\;(\partial^*X,(d)^p_\varepsilon)\to (\partial^*Y,(d')^{p'}_\varepsilon)$$ is $\eta$-quasisymmetric with $\eta=\eta_{\delta, \lambda, \mu}$.
\end{Thm}

A space is {\it a length space} if the distance between any two points in this space is equal to the infimum of the lengths of all curves joining these two points.
A mapping $f:$ $(X,d)\to (Y,d')$ is {\it weakly surjective} if for any fixed $q\in Y$, $$\limsup_{d'(y,q)\to \infty}\frac{d'(y, f(X))}{d'(y,q)}<1.$$ Obviously,
surjectivity implies weak surjectivity.

\subsection{Notational conventions}
If there is no danger of confusion, the metric $d$ will be dropped from all involved notations. For example, we write $k(d)_X=k_X$, $\ell_d=\ell$, $\delta(d)_X=\delta_X$, $\sigma(d)_X=\sigma_X$,  and so on.
In particular, $G$ always denotes a proper subdomain of $\IR^n$.
Since $k(\sigma(d))_G=k(d)_G$ when $d=|\cdot|$, the Euclidean metric, we simply use $k_G$ to denote both of them. Also, we simply use $\sigma$ to denote $\sigma_G$.

\subsection{Rough starlikeness}
Let $(X,d)$ be a length Gromov $\delta$-hyperbolic space, and let $\mu$ and  $h$ be nonnegative constants.
A {\it $(\mu,h)$-road} $\overline{\alpha}$ in $X$ is a sequence of arcs $\alpha_i$ with endpoints $y_i$ and $u_i$ along the direction from $y_i$ to $u_i$ satisfying the following:
\ben
\item
each $\alpha_i$ is $h$-short;
\item
the sequence of lengths $\ell_d(\alpha_i)$ is increasing and tending to $\infty$;
\item
for $i\leq j$, the length mapping (cf. \cite[Subsection 2.16]{Vai10}) $g_{ij}: \alpha_i\to \alpha_j$ with $g_{ij}(y_i)=y_j$ satisfies $d(g_{ij}(x), x)\leq \mu$ for all $x\in \alpha_i$.
\een

For $x$, $y\in X$, and $h\geq 0$, a curve $\gamma$ in $X$ connecting $x$ and $y$ is called {\it $h$-short} if $$\ell(\gamma)\leq d(x,y)+h.$$

By \cite[Lemma 6.3]{Vai10}, we see that for a $(\mu,h)$-road $\overline{\alpha}$ which consists of the arcs $\alpha_i$ connecting $y_i$ and $u_i$ along the direction from $y_i$ to $u_i$, the corresponding sequence $\{u_i\}$ is Gromov and defines a point $b$ on $\partial^*X$. If, further, for each $i$, $y_i=y$, then we say that $\overline{\alpha}$ is a road connecting $y$ and $b$.

Now, we are ready to state two definitions for rough starlikeness. The first one is as follows.
Let $(X, d)$ be a Gromov $\delta$-hyperbolic space, and let $K$, $\mu$ and $h$ be nonnegative constants. We say that $X$ is
 \ben
 \item
 {\it $(K,\mu, h)$-roughly starlike} with respect to a base point $w\in X$ if for any $x\in X$, there is a $(\mu,h)$-road $\overline{\alpha}$ connecting $w$ and
some point $b\in \partial^*X$ such that $d(x,\overline{\alpha})\leq K$.
 \item
{\it $(K,\mu)$-roughly starlike} with respect to a base point $w\in X$ if it is $(K,\mu, h)$-roughly starlike with respect to $w$ for all $h>0$.
 \een

The following is another definition for rough starlikeness.

A proper, geodesic and Gromov hyperbolic space $(X, d)$ is said to be {\it $K$-roughly starlike} $(K\geq 0)$ with respect to a base point $w\in X$ if for each point $x\in X$, there exists a geodesic ray $\beta$ starting from $w$ such that $d(x,\beta)\leq K$.

We shall prove that these two definitions for rough starlikeness are equivalent when the spaces are proper, geodesic and Gromov hyperbolic (see Lemma \ref{tue-1} below). First, let us recall a result due to V\"ais\"al\"a.

\begin{Thm}\label{Lem-03} $($\cite[Theorem $6.32$]{Vai10}$)$
Suppose that $(X, d)$ is a length Gromov $\delta$-hyperbolic space. Let $\varphi:[0,\infty)\to X$ be a $(\lambda,\mu)$-roughly quasi-isometric curve, and let $\overline{\alpha}$ be a $(\mu,h)$-road connecting $\varphi(0)$ and $\varphi(\infty)$. Then $$d_H(\overline{\alpha},\;\varphi)\leq M,$$ where $d_H$ stands for the Hausdorff distance and $M=M(\delta,\lambda,\mu,h)$.
\end{Thm}

\blem\label{tue-1}
Suppose that $(X, d)$ is a proper, geodesic and Gromov $\delta$-hyperbolic space. Then the following are equivalent:
\ben
\item\label{tue-1(1)}
$X$ is $(K_1, \mu_1, h_1)$-roughly starlike;
\item\label{tue-1(2)}
 $X$ is $(K_2, \mu_2)$-roughly starlike;
 \item\label{tue-1(3)}
 $X$ is $K_3$-roughly starlike.
 \een
\elem
\bpf Since \cite[Lemma $6.34(1)$]{Vai10} implies that the conditions \eqref{tue-1(1)} and \eqref{tue-1(2)} in the lemma are equivalent, and since the implication from \eqref{tue-1(3)} to \eqref{tue-1(1)} is obvious, we see that, to prove this lemma, it suffices to show the implication from \eqref{tue-1(1)} to \eqref{tue-1(3)}. To this end, let $w\in X$. Then the assumption guarantees that there is a $(\mu_1, h_1)$-road $\overline{\alpha}$ connecting $w$ and $b\in \partial^*X$ with $d(w,\overline{\alpha})\leq K_1$. Furthermore, since $(X, d)$ is proper, geodesic and Gromov $\delta$-hyperbolic, we see that Hopf-Rinow Theorem is applicable in this situation (see, e.g., \cite[Lemma 3.1 in Part III-H]{BH99}). This theorem implies that there is a geodesic ray $\beta$ connecting $w$ and $b$.
 Thus, by Theorem \Ref{Lem-03}, we know that there is a constant $M=M(\delta, \mu_1, h_1)$ such that $d_H(\overline{\alpha}, \beta)\leq M$, which leads to $$d(w,\beta)\leq M+K_1.$$ Now, the lemma follows by letting $K_3=M+K_1$.
\epf

Let us recall the following results concerning rough starlikeness, which are useful for the discussions in Section \ref{sec-3}.

\begin{Thm}\label{fri-1} $($\cite[Theorem 3.22]{Vai05}$)$
Every Gromov $\delta$-hyperbolic domain in Banach spaces is $(K,\mu)$-roughly starlike with
respect to each point in this domain, where $K=K(\delta)$ and $\mu=4\delta+1$.\end{Thm}

%We remark from Theorem \Ref{fri-1} it follows that every $\delta$-hyperbolic domain $X$ in a Euclidean space is roughly starlike
%with constant depending only on $\delta$. Indeed, in this case every hyperbolic domain is a proper complete geodesic space; see \cite[$2.8$]{BHK}.

\begin{Thm}\label{Lem-4}$($\cite[Proposition $4.37$]{BHK}$)$
If $(X,d)$ is a $K$-roughly starlike, proper, geodesic and Gromov $\delta$-hyperbolic space, then for any $0<\varepsilon\leq \varepsilon_0$, the identity mapping
from $(X,d)$ to $(X,k((d)_{\varepsilon,p})_X)$ is homeomorphic and $M_0$-quasi-isometric, where $\varepsilon_0$ and $p\in X$ are from Theorem \Ref{Lem-1} and $M_0=M_0(K,\delta,\varepsilon)$.
\end{Thm}

\begin{Thm}\label{Lem-3}$($\cite[Theorem $3.6$]{BHK}$)$
If $(X, d)$ is a uniform space, then $(X, k(d)_X)$ is a proper, geodesic and Gromov $\delta$-hyperbolic space. Moreover, if $X$ is bounded, then $(X, k(d)_X)$ is roughly starlike, and the quasisymmetric gauge determined by $d$ on the metric boundary $\partial X$ is naturally equivalent to the canonical gauge on the Gromov boundary $\partial^*X$.
\end{Thm}

Note that the canonical gauge in Theorem \Ref{Lem-3} consists of visual metrics on $\partial^*G$. See \cite{BHK} for the details.

\subsection{Two conditions}\label{prop-AB}

Let us start with the introduction of several notations.
Assume that $(G,|\cdot|)$ is a bounded domain in $\mathbb{R}^n$. Let $w_0\in G$ be such that
\be\label{wen-2} \delta_G(w_0)=\max\{\delta_G(x):\; x\in G\},\ee
and let
\be\label{shantou-4} \nu_0=\min\{1,\varepsilon_0, \;\tau_0,\; \mu_0\},\ee
where $\varepsilon_0$ (resp. $\tau_0$ and $\mu_0$) is defined by Theorem \Ref{Lem-1} (resp. Subsection \ref{changsha-1} and Theorem \Ref{sat-12}).

Next, we introduce the following two conditions.

\bcond\label{PropertyA}
We say that $G$, a bounded domain in $\mathbb{R}^n$ $(n\geq 2)$, satisfies {\it Condition \ref{PropertyA}} if
\ben\item\label{PropertyA-1}
$(G, k_G)$ is Gromov $\delta$-hyperbolic, and
\item\label{PropertyA-2}
if $w_0\in G$ is as in \eqref{wen-2}, then there exists a bijective natural mapping
$$\varphi:\; (G^*, (k_G)_{\tau}^{w_0})\to (\overline{G_{\sigma}}, \sigma)$$
such that the restriction $\varphi|_{\partial_{\sigma} G}:$ $(\partial^* G, (k_G)_{\tau}^{w_0})\to (\partial_{\sigma} G, \sigma)$ is $\eta$-quasisymmetric, where $\tau\in(0,\nu_0)$.
\een
\econd

\bcond\label{PropertyB}
We say that $G$, a bounded domain in $\mathbb{R}^n$ $(n\geq 2)$, satisfies {\it Condition \ref{PropertyB}} if there is a homeomorphism $f:$ $(\overline{G_{\sigma}},\sigma)\to (\overline{G'},d')$, where $G'$ is an $M$-uniform domain in $(Y,d')$, a metric space, such that
\ben
\item
the restriction $f|_{G}:$ $(G,\sigma)\to (G',d')$ is $M$-quasihyperbolic, and
\item\label{PropertyB-3}
the restriction $f|_{\partial_{\sigma}{G}}:$ $(\partial_{\sigma}{G},\sigma)\to (\partial_{d'}{G'},d')$ is $\eta$-quasisymmetric.
\een
\econd

Obviously, Condition \ref{PropertyA} coincides with the condition $(\mathfrak{ii})$ in Theorem \ref{thm-1.1}. Then the equivalence of the conditions $(\mathfrak{i})$ and $(\mathfrak{ii})$ in Theorem \ref{thm-1.1} easily follows from the following theorem.

\begin{thm}\label{thm-2.3}
Suppose that $G$ is a bounded domain in $\mathbb{R}^n$ $(n\geq 2)$. Then the following statements are equivalent:\ben
\item[$(\mathfrak{a})$]\label{thu-4}
$G$ is inner uniform;
\item[$(\mathfrak{b})$]\label{thu-5}
$G$ satisfies Condition \ref{PropertyA};
\item[$(\mathfrak{c})$]\label{thu-6}
$G$ satisfies Condition \ref{PropertyB}.
\een
\end{thm}

\medskip

The implications of $(\mathfrak{a})$ $\Longrightarrow$ $(\mathfrak{b})$ $\Longrightarrow$ $(\mathfrak{c})$  will be proved in Section \ref{sec-3} (see Proposition \ref{thm-2.1}), and the implication $(\mathfrak{c})$ $\Longrightarrow$ $(\mathfrak{a})$ will be shown in Section \ref{sec-6} (see Theorem \ref{thm-6.1}).

%%%%%%%%%%%%%%%%%%%%%%%%%%%%%%%%%%%%%%%%%%%
\section{Proofs of the implications $(\mathfrak{a})$ $\Longrightarrow$ $(\mathfrak{b})$ $\Longrightarrow$ $(\mathfrak{c})$ in Theorem \ref{thm-2.3}}\label{sec-3}
%%%%%%%%%%%%%%%%%%%%%%%%%%%%%%%%%%%%%%%%%%%

The purpose of this section is to prove the following result.

\bprop\label{thm-2.1} Suppose that $G$ is a bounded domain in $\mathbb{R}^n$ $(n\geq 2)$. Then the implications $(\mathfrak{a})$ $\Longrightarrow$ $(\mathfrak{b})$ $\Longrightarrow$ $(\mathfrak{c})$ in Theorem \ref{thm-2.3} are true.
\eprop

\subsection*{Proof of the implication $(\mathfrak{a})$ $\Longrightarrow$ $(\mathfrak{b})$}
Assume that $G$ is inner uniform. Then the first assertion in \cite[Theorem 1.11]{BHK} implies that the first statement of Condition \ref{PropertyA} is true. Since Lemma \ref{changsha-2} ensures that $G$ is also bounded with respect to the inner metric, by replacing the point $p$ by the one $w_0$ defined in \eqref{wen-2}, we know from Theorem \Ref{sat-12} that the second statement of Condition \ref{PropertyA} holds true as well.

\subsection*{Proof of the implication $(\mathfrak{b})$ $\Longrightarrow$ $(\mathfrak{c})$}
Assume that $G$ satisfies Condition \ref{PropertyA}. To prove this implication, it is enough to establish the following theorem.

\bthm\label{lem-1}
There exists a homeomorphism $f:$ $(\overline{G_{\sigma}}, \sigma)\to \big(\overline{G_{(k_G)_{\tau,w_0}}}, (k_G)_{\tau,w_0}\big)$, where $\tau\in(0,\nu_0)$, such that
 \ben
 \item[$(i)$]
 $(G,(k_G)_{\tau,w_0})$ is an $M$-uniform space;
 \item[$(ii)$]
  the restriction $f|_{G}=\mathop\mathrm{id}$: $(G,\sigma)\to (G,(k_G)_{\tau,w_0})$ is $M$-quasihyperbolic, and
 \item[$(iii)$]
the restriction $f|_{\partial_{\sigma} G}$: $(\partial_{\sigma} G, \sigma)\to (\partial_{(k_G)_{\tau,w_0}} G, (k_G)_{\tau,w_0})$ is $\eta$-quasisymmetric.
 \een\ethm

We make the following notational conventions: In Theorem \ref{lem-1} and its proof below, we always use the symbol $\eta$ (resp. the symbols $\delta$ and $K$) to denote the control functions of quasisymmetry (resp. the coefficients of Gromov hyperbolicity and rough starlikeness). We also use $M$ to denote the coefficients of uniformity or quasi-isometry or quasihyperbolicity. These functions and constants depend only on the given parameters. Also, they are not necessarily the same on different occasions.

The following lemma is useful.

\blem\label{wen-1}
The space $(G, k_G)$ is $K$-roughly starlike, complete, proper, geodesic and Gromov $\delta$-hyperbolic.\elem
\bpf  Obviously, the identity mapping $(G,|\cdot|)\to (G,\sigma)$ is a local isometric homeomorphism. We see from \cite[Proposition 2.8]{BHK} that $(G, k_G)$ is complete,  proper, geodesic and Gromov $\delta$-hyperbolic.
Moreover, Lemma \ref{tue-1} and Theorem \Ref{fri-1} guarantee the rough starlikeness of $(G, k_G)$, and thus, the lemma is proved.\epf

\subsection*{Proof of Theorem 3.1}

\bpf
First, we know from Lemma \ref{wen-1} and Theorem \Ref{Lem-1} that:

$(\mathfrak{S}_1)$\label{mon-22}
the space $(G,(k_G)_{\tau,w_0})$ is bounded and $M$-uniform.

As a direct consequence of the statement $(\mathfrak{S}_1)$ and Theorem \Ref{Lem-3}, we have that

$(\mathfrak{S}_2)$\label{thu-1}
the space $\big(G, k\big((k_G)_{\tau,w_0}\big)_G\big)$ is $K$-roughly starlike, proper, geodesic and Gromov $\delta$-hyperbolic.

Now, we are ready to start constructing the required homeomorphism. First, it follows from Lemma \ref{wen-1} and Theorem \Ref{Lem-4} that the following identity mapping
$$\mathop{\mathrm{id}_1}:\;(G, k_G)\to \big(G, k\big((k_G)_{\tau,w_0}\big)_G\big)$$
is homeomorphic and $M$-quasi-isometric. Hence, by Theorem \Ref{fri-2}, together with Lemma \ref{wen-1} and the statement $(\mathfrak{S}_2)$, we see that there exists a bijective natural mapping
$$\mathop{\mathrm{id}_1^*}: \;(G^*, (k_G)_{\tau}^{w_0})\to \Big(G^*, \Big(k\big((k_G)_{\tau,w_0}\big)_G\Big)_{\tau}^{w_0}\Big)$$
such that the restriction
$$\mathop{\mathrm{id}_1^*}|_{\partial^* G}:\;(\partial^* G, (k_G)_{\tau}^{w_0})\to \Big(\partial^* G, \Big(k\big((k_G)_{\tau,w_0}\big)_G\Big)_{\tau}^{w_0}\Big)$$
is $\eta$-quasisymmetric.

Moreover, we infer from the statement $(\mathfrak{S}_1)$ and Theorem \Ref{qz-2} that there is a bijective positive mapping
 \be\label{qz-6}
 \psi:\; \Big(G^*,\Big(k\big((k_G)_{\tau,w_0}\big)_G\Big)_{\tau}^{w_0}\Big)\to \big(\overline{G_{(k_G)_{\tau,w_0}}}, (k_G)_{\tau,w_0}\big).\ee

In order to exploit Theorem \Ref{sat-12}, we note the following relationships between the distance from $w_0$ to $\partial_{(k_G)_{\tau,w_0}} G$ and the diameter of $G$ with respect to $(k_G)_{\tau,w_0}$:
\beqq
\delta((k_G)_{\tau,w_0})_G(w_0) &\geq &  \frac{1}{\tau e}\rho(k_G)_{\tau,w_0}(w_0)  \geq  \frac{1}{2e}\diam_{(k_G)_{\tau,w_0}} (G),
\eeqq
where the first inequality follows from \cite[$(4.6)$]{BHK}, and the second one is from the inequality next to \cite[$(4.3)$]{BHK}.
 Hence
for all $x\in G$,
\beqq\label{qz-3}\delta((k_G)_{\tau,w_0})_G(x)\leq 2e\,\delta((k_G)_{\tau,w_0})_G(w_0),\eeqq
and thus, we see from Theorem \Ref{sat-12} that
 the restriction
 $$\psi|_{\partial^* G}:\; \Big(\partial^* G,\Big(k\big((k_G)_{\tau,w_0}\big)_G\Big)_{\tau}^{w_0}\Big)\to (\partial_{(k_G)_{\tau,w_0}} G, (k_G)_{\tau,w_0})$$
is $\eta$-quasisymmetric.

Let
$$f_1=\psi\circ \mathop{\mathrm{id}_1^*}\circ \varphi^{-1}|_{\partial_{\sigma} G}:\;(\partial_{\sigma} G, \sigma)\to (\partial_{(k_G)_{\tau,w_0}} G, (k_G)_{\tau,w_0}),$$ where the mapping $\varphi$ is from Condition \ref{PropertyA}.
Then $f_1$ is $\eta$-quasisymmetric
since a composition of quasisymmetric mappings is still quasisymmetric.

Furthermore, let
$\mathop{\mathrm{id}_2}:$ $(G, \sigma)\to (G, k_G)$ and $\mathop{\mathrm{id}_3}:$ $(G, (k_G)_{\tau,w_0})\to (G, k((k_G)_{\tau,w_0})_G)$ be identity mappings.
By the statement $(\mathfrak{S}_1)$, we see that $(G,(k_G)_{\tau,w_0})$ is uniform. Then it follows from \cite[Proposition 2.8]{BHK} that both $\mathop{\mathrm{id}_2}$ and $\mathop{\mathrm{id}_3}$ are homeomorphisms. Let
$$f_2=\mathop{\mathrm{id}_3^{-1}}\circ \mathop{\mathrm{id}_1}\circ \mathop{\mathrm{id}_2}:\; (G,\sigma)\to (G, (k_G)_{\tau,w_0}).$$
Then $f_2$ is also a homeomorphism. Further, we know that $f_2$ is $M$-quasihyperbolic since $\mathop{\mathrm{id}_1}$ is $M$-quasi-isometric.
Let
\beqq\label{hoemo} f:\; (\overline{G_{\sigma}}, \sigma)\to \big(\overline{G_{(k_G)_{\tau,w_0}}}, (k_G)_{\tau,w_0}\big)\eeqq be defined as follows:
$$f|_G=f_2\;\;\mbox{and}\;\;f|_{\partial_{\sigma} G}=f_1.$$

Next, we prove that the mapping $f$ is a homeomorphism. Because $f$ is bijective, it is sufficient to verify the continuity of both $f$ and $f^{-1}$. For this, we need the following assertion.

$(\mathfrak{S}_3)$\label{z1} Suppose that $\overline{x}=\{x_n\}\subset G$, $\lambda\in\partial^* G$ and $\xi\in\partial_\sigma G$ with $\varphi(\lambda)=\xi$. Then $\sigma(x_n,\xi)\to 0$ as $n\to \infty$, if and only if $\overline{x}$ is a Gromov sequence in $(G,k_G)$ and $\overline{x}\in \lambda$.

The sufficiency part in the assertion $(\mathfrak{S}_3)$ follows from the assumption that $G$ satisfies Condition \ref{PropertyA} and from Theorem \Ref{qz-1}. To check the necessity part in $(\mathfrak{S}_3)$, assume that $\sigma(x_n,\xi)\to 0$ as $n\to \infty$. Then we claim that $\overline{x}$ is a Gromov sequence in $(G,k_G)$. To prove this assertion, let $\alpha_{n,m}$ denote a quasihyperbolic geodesic in $G$ connecting $x_n$ and $x_m$.

First, since $G$ is assumed to be Gromov $\delta$-hyperbolic, we know from Theorem \Ref{BB03} that it satisfies the Gehring-Hayman condition. This means that there exists a constant $c_{gh}\geq 1$ such that for any positive integers $n$ and $m$,
\be\label{21-10-6}
\ell(\alpha_{n,m})\leq c_{gh}\sigma(x_n,x_m),
\ee where $c_{gh}$ is independent of the subscripts $n$ and $m$.

Furthermore, by Lemma \ref{wen-1} and \cite[$(3.2)$]{BHK}, we see that there is a point $z_{n,m}\in\alpha_{n,m}$ such that
\be\label{21-10-7}
(x_n|x_m)_{w_0}\geq k_G(z_{n,m},w_0)-8\delta\geq \Big|\log \frac{\delta(\sigma)_G(w_0)}{\delta(\sigma)_G(z_{n,m})}\Big|-8\delta.
\ee

Based on the fact
$$\delta(\sigma)_G(z_{n,m})\leq \sigma(z_{n,m},\xi)\leq \ell(\alpha_{n,m})+\sigma(x_n,\xi),$$
we see from \eqref{21-10-6} that
$$\delta(\sigma)_G(z_{n,m})\to 0$$ as $n, m\to\infty$, and thus, \eqref{21-10-7} implies
$$(x_n|x_m)_{w_0}\to\infty$$ as $n, m\to\infty$. This ensures that $\overline{x}$ is a Gromov sequence, and hence, the assertion is proved.

Since $\overline{x}$ is a Gromov sequence, we know that there is a $\zeta\in\partial^* G$ such that $\overline{x}\in \zeta$.
As $\varphi$ is bijective, we deduce from Theorem \Ref{qz-1} that $\zeta=\lambda$, which gives the sufficiency part in $(\mathfrak{S}_3)$.\medskip

Now, let us return to the proof of the theorem. Let $\xi\in\partial_\sigma G$ and $\overline{x}=\{x_n\}\subset G$. By the assertion $(\mathfrak{S}_3)$, we know that $\sigma(x_n,\xi)\to 0$ as $n\to \infty$, if and only if $\overline{x}$ is a Gromov sequence in $(G,k_G)$ and $\overline{x}\in \lambda_1$, where $\lambda_1=\varphi^{-1}(\xi)$.  Moreover, by Lemma \ref{wen-1} 
and the statement $(\mathfrak{S}_2)$ along with
 the fact that $\mathop{\mathrm{id}_1}$ is homeomorphic and $M$-quasi-isometric, we see from \cite[Proposition $6.3$]{BS} that $\overline{x}$ is a Gromov sequence in $(G,k_G)$ and $\overline{x}\in \lambda_1$, if and only if $\overline{x}$ is a Gromov sequence in $(G,k((k_G)_{\tau,w_0})_G)$ and $\overline{x}\in \lambda_2$, where $\lambda_2=\mathop{\mathrm{id}_1^*}(\lambda_1)$.

 Furthermore, we deduce from the statement $(\mathfrak{S}_1)$ and Theorem \Ref{qz-2} that $\overline{x}$ is a Gromov sequence in $(G,k((k_G)_{\tau,w_0})_G)$ and $\overline{x}\in \lambda_2$, if and only if $(k_G)_{\tau,w_0}(x_n,\xi_1)\to 0$, where $\xi_1=\psi(\lambda_2)$ and $\psi$ is from \eqref{qz-6}. Since $\xi_1=f(\xi)$, we see that both $f$ and $f^{-1}$ are continuous.
Hence the proof of the theorem is complete.
\epf

%%%%%%%%%%%%%%%%%%%%%%%%%%%%%%%%%%%%%%%%%%%
%%%%%%%%%%%%%%%%%%%%%%%%%%%%%%%%%%%%%%%%%%%
\section{Auxiliary lemmas} \label{sec-4}
%%%%%%%%%%%%%%%%%%%%%%%%%%%%%%%%%%%%%%%%%%%
%%%%%%%%%%%%%%%%%%%%%%%%%%%%%%%%%%%%%%%%%%%

This section consists of three parts. In the first part, we give explicit expressions of certain constants which will be used in the discussions in Sections \ref{sec-4}$-$\ref{sec-6} of the paper. In the second part,
we recall a series of results from \cite{HRWZ} which will be applied later on. In the last part, we introduce a new condition, and prove a related lemma.

%%%%%%%%%%%%%%%%%%%%%%%%%%%%%%%%%%%%%%%%%%%

\subsection{Several constants}  \label{sec-4-0}

In the rest of this paper, we assume that $G$ is a bounded domain in $\IR^n$ $(n\geq 2)$,   $f:$ $(\overline{G_{\sigma}},\sigma)\to (\overline{G'},d')$ a homeomorphism, $f|_{\partial_{\sigma}{G}}:$ $(\partial_{\sigma}{G},\sigma)\to (\partial_{d'}{G'},d')$ is an $\eta$-quasisymmetric mapping with $\eta(1)\geq 1$,  $f|_{G}:$ $(G,\sigma)\to (G',d')$ is an $M$-quasihyperbolic mapping, $G'$ is an $M$-uniform domain with $M\geq 37$  and $\lambda_1$ is the coefficient of the $\lambda_1$-quasigeodesic appeared in the following section with $1\leq\lambda_1\leq 100M^2$.

For convenience, we list several commonly used constants in the rest of the paper.
\ben
\item\label{HWRZ-25-0}  $\nu_1=\nu_1(\lambda_1,M)=e^{(2\lambda_1M)^{10}}$ with $\lambda_1\geq 1$, where
$\nu_1$ is from Lemma \Ref{h-W-1};

\item\label{HWRZ-2}  $\nu_2=M+\nu_1(\lambda_1,M)$ with $\lambda_1=100 M^2$, where
$\nu_2$ is from Lemmas \Ref{lab-40} and \Ref{shantou-5}; 

\item\label{HWRZ-3} $\nu_3=e^{(2\nu_2M)^2}$, where $\nu_3$ is from Condition C;

\item\label{HWRZ-25-1} $M_2=10\nu_3^4\eta(\nu_3)\max\Big\{1,\big(\eta^{-1}(\frac{1}{5}\nu_2^{-3})\big)^{-1}\Big\}$, where $M_2$ is from Lemma \Ref{H-W-Zhou-1};

\item\label{HWRZ-4} $M_1=\max\Big\{\eta(M_2^5)\cdot e^{MM_2^5},\big(\eta^{-1}(M_2^{-1})\big)^{-1}M_2\Big\}$, where $M_1$ is from Condition B;

\item\label{HWRZ-5} $M_0=\eta(\nu_3M_1^2)\cdot e^{\nu_3M_1}$, where $M_0$ is from the proof
of Lemma \ref{Zhou-W-0};

\item\label{HWRZ-6} $B_0=20M_0^2$, where $B_0$ is from Theorem \Ref{Thm-HRWZ}.
\een
\medskip

\br\label{nw-1}
By the list as above, the following useful relations are obvious: \ben \item\label{HWRZ-25-2} $\nu_2 \geq e^{(100M^3)^{10}}$ and $\nu_3\geq e^{(74\nu_2)^2}$;

\item\label{HWRZ-25-3} $M_2\geq10\nu_3^4 $  and $M_1>e^{37M_2^5}\max\{1, \eta(e^{10\nu_3^4})\}$;

\item\label{HWRZ-25-4} $M_0\geq \max\{e^{\nu_3M_1}, \eta(e^{10\nu_3^3})e^{\nu_3M_1}\}$ and $B_0\geq \max\{M_0^\frac{3}{2}(4M_1)^6, (2M_1)^{10}\}$.
\een
\er

%%%%%%%%%%%%%%%%%%%%%%%%%%%%%%%%%%%%%%%%%%%

\subsection{Known results} In this subsection, we recall a series of results from \cite{HRWZ}.

%\begin{Lem}\label{H-W-Zhou-0} $($\cite[Lemma 3.1]{HRWZ}$)$ Suppose that $x_1$ and $x_2$ are two points in $(X,d)$.
%Let $\alpha_{12}$ denote a curve in $X$ connecting
%$x_1$ and $x_2$.  If there is a constant $a\geq 1$ such that for any $x\in\alpha_{12}$,
%$\ell(\alpha_{12}[x_1,x])\leq a \delta_X(x)$,
%then $$k_X(x_1,x_2)\leq 4a\log\Big(1+\frac{\ell(\alpha_{12})}{\delta_X(x_1)}\Big).$$
%\end{Lem}

%%%%%%%%%%%%%%%%%%%%%%%%%%%%%%%%%%%%%%%%%%%%

\begin{Lem}\label{h-W-1}$($\cite[Lemma 3.2]{HRWZ}$)$ Suppose that $(X,d)$ is $M$-uniform.
\ben
\item\label{h-W-1(1)}
For all $u,$ $v\in X$,
$$k_X(u,v)\leq 4M^2\log\Big(1+\frac{d(u,v)}{\min\{\delta_X(u),\delta_X(v)\}}\Big).$$
\item\label{h-W-1(2)}
Every $\lambda_1$-quasigeodesic or $\lambda_1$-quasigeodesic ray $\alpha$ in $X$ with $\lambda_1\geq 1$ is $\nu_1$-uniform, where $\nu_1=\nu_1(\lambda_1,M)=e^{(2\lambda_1M)^{10}}$.
\een
\end{Lem}

\begin{Lem}\label{lab-40}$($\cite[Lemma 3.3]{HRWZ}$)$
Suppose that $(X,d)$ is $M$-uniform and $\alpha$ is a quasihyperbolic geodesic in $X$. Then
for all $u$, $v\in \alpha$, $\alpha[u,v]$ is $\nu_2$-uniform and $\ell(\alpha[u,v])\leq \nu_2d(u,v)$, where $\nu_2=M+\nu_1(100M^2,M)$.
\end{Lem}

%%%%%%%%%%%%%%%%%%%%%%%%%%%%%%%%%%%%%%%%%%%%

\begin{Lem}\label{NY-1}$($\cite[Lemma 3.4]{HRWZ}$)$ Suppose that $f:$ $(X,d)\to (Y,d')$ is $M$-quasihyperbolic, and $\alpha$ is a $\lambda_1$-quasigeodesic $($resp. a $\lambda_1$-quasigeodesic ray$)$ in $X$ with $\lambda_1\geq 1$. Then
 $\alpha'=f(\alpha)$ is a $\lambda_1M^2$-quasigeodesic $($resp. a $\lambda_1M^2$-quasigeodesic ray$)$.
 \end{Lem}

\begin{Lem}\label{shantou-5}$($\cite[Lemma 3.5]{HRWZ}$)$
Suppose that $f:$ $(X,d)\to (Y,d')$ is $M$-quasihyperbolic, $(Y,d')$ is $M$-uniform and $\alpha$ is a $\lambda_1$-quasigeodesic ray in $X$ with $1\leq \lambda_1\leq 100$. Then $\alpha'=f(\alpha)$ is $\nu_2$-uniform,  where $\nu_2=M+\nu_1(100M^2,M)$.
\end{Lem}

%%%%%%%%%%%%%%%%%%%%%%%%%%%%%%%%%%%%%%%%%%%%

\begin{Cons}\label{Conditions-A} $($\cite[Conditions A]{HRWZ}$)$
Suppose that $x_1,$ $x_2$ and $x_3$ are points in $(G,|\cdot|)$, and
 $\alpha_{23}$ denotes a curve in $G$ with endpoints $x_2$ and $x_3$.
  We say that the quadruple $[x_1,x_2,x_3;\alpha_{23}]$ satisfies {\it Condition A} if
\ben
\item\label{H-W-Zh-X-1(1)}
 $\sigma(x_2,x_3)\geq 20\sigma(x_1,x_2)$, and
\item\label{H-W-Zh-X-1(2)}
 $k_{G}(x_1,\alpha_{23})\geq 10$, where $k_{G}(x_1,\alpha_{23})$ denotes the quasihyperbolic distance from $x_1$ to the curve $\alpha_{23}$.
\een
\end{Cons}

\noindent Recall that $G$ is a proper subdomain of $\IR^n$ and $|\cdot|$ denotes the Euclidean metric.
%%%%%%%%%%%%%%%%%%%%%%%%%%%%%%%%%%%%%%%%%%%

\begin{Lem}\label{H-W-Zh-X-1} $($\cite[Lemma 3.6]{HRWZ}$)$
 If the quadruple $[x_1,x_2,x_3;\alpha_{23}]$ satisfies Condition A, then there exists a $100$-quasigeodesic ray $\alpha$ in $G$ starting from $x_1$
and ending at $x_{1,1}\in \partial_{\sigma}G$ such that for any $x\in\alpha$,
$$\frac{1}{42}\sigma(x_1,x_2)< \sigma(x_2,x)\leq 5\sigma(x_1,x_2)\;\;\mbox{and}\;\;\ell(\alpha[x,x_{1,1}))\leq 5\delta_{G}(x),$$
where $\alpha[x,x_{1,1})$ denotes the sub-ray of $\alpha$, which starts from $x$ and ends at $x_{1,1}$.
\end{Lem}

%%%%%%%%%%%%%%%%%%%%%%%%%%%%%%%%%%%%%%%%%%%%

\begin{Cons}\label{Conditions-B} $($\cite[Conditions B]{HRWZ}$)$
Suppose that $(X,d)$ is $M$-uniform. Let $x_1$, $x_2\in X$, $x_3\in \overline{X_d}$ and
 $\alpha_{23}$ a curve in $X$ with endpoints $x_2$ and $x_3$ or a ray in $X$ starting from $x_2$ and ending at $x_3\in \partial_d X$. We say that the quadruple $[x_1,x_2,x_3;\alpha_{23}]$
 satisfies {\it Condition B} if \ben
  \item\label{changsha-4}
 $d(x_1,x_2)\leq 2\nu_2^2\delta_{X}(x_1)$,
 \item\label{changsha-5}
 $k_{X}(x_1, \alpha_{23})> \frac{1}{20}\nu_2^{-1}\log M_1,$ and
 \item\label{changsha-6}
 $\alpha_{23}$ is $\nu_2$-uniform.
  \een
 \end{Cons}
%%%%%%%%%%%%%%%%%%%%%%%%%%%%%%%%%%%%%%%%%%%%

\begin{Lem}\label{H-W-Zhou-1} $($\cite[Lemma 3.9]{HRWZ}$)$    If the quadruple $[x_1,x_2,x_3;\alpha_{23}]$ satisfies Condition B,
then
$$\delta_{X}(x_1)> M_2\ell(\alpha_{23}).$$
\end{Lem}

%%%%%%%%%%%%%%%%%%%%%%%%%%%%%%%%%%%%%%%%%%%%

%\noindent{\bf Property $F$.}
\begin{Cons}\label{Conditions-C}$($\cite[Conditions C]{HRWZ}$)$
Suppose that $(X,d)$ is $M$-uniform.
Let $x_1$,  $x_2$ and $x_3\in X,$ $x_4\in \overline{X_d}$, and let $\alpha_{12}$ be a curve in $X$ with endpoints $x_1$ and $x_2$, $\alpha_{24}$ a curve in $X$ with endpoints $x_2$ and $x_4$ or a ray in $X$ starting from $x_2$ and ending at $x_4\in \partial_d X$. We say that the sextuple $[x_1,x_2,x_3,x_4;\alpha_{12},\alpha_{24}]$ satisfies {\it Condition C} if \ben
\item\label{H-W-2-zhou(1)}
$d(x_1,x_4)\leq \frac{1}{2}\nu_3^{-2}\min\{d(x_1,x_2),d(x_2,x_4)\}$,
\item\label{H-W-2-zhou(3)}
$\alpha_{12}$ is a quasihyperbolic geodesic such that for any $x\in\alpha_{12}$, $\ell(\alpha_{12}[x_1,x])\leq \nu_2\delta_{X}(x),$
\item\label{H-W-2-zhou(4)}
$x_3\in \alpha_{12}$ and $\ell(\alpha_{12}[x_1,x_3])\geq 2d(x_1,x_4)$, and
\item\label{H-W-2-zhou(2)}
$\alpha_{24}$ is $\nu_2$-uniform.
\een
\end{Cons}
%%%%%%%%%%%%%%%%%%%%%%%%%%%%%%%%%%%%%%%%%%%%

\begin{Lem}\label{H-W-2-zhou} $($\cite[Lemma 3.10]{HRWZ}$)$ Suppose that the sextuple $[x_1,x_2,x_3,x_4;\alpha_{12},\alpha_{24}]$ satisfies Condition C. Then $$k_{X}(x_3, \alpha_{24})\leq 12M\log \nu_3.$$
\end{Lem}

%%%%%%%%%%%%%%%%%%%%%%%%%%%%%%%%%%%%%%%%%%%

The following is one of the main results in \cite{HRWZ}. It is useful for the arguments in Section \ref{sec-6}.

\begin{Thm}\label{Thm-HRWZ}$($\cite[Theorem 4.1]{HRWZ}$)$  Suppose that $G$ satisfies Condition \ref{PropertyB}. Let $\gamma\subset G$ be a curve connecting two points $u$ and $v$ in G. If the image $\gamma'=f(\gamma)$ is a quasihyperbolic geodesic in $G'$, then for any $x\in\gamma$,
$$\min\{\diam(\gamma[u,x]), \diam(\gamma[v,x])\}\leq B_0\delta_G(x),$$
where $B_0$ is a constant depending only on the given parameters $($See Subsection \ref{sec-4-0}$)$.
\end{Thm}

%%%%%%%%%%%%%%%%%%%%%%%%%%%%%%%%%%%%%%%%

\subsection{A new condition}
In this subsection, we introduce a new condition concerning points and curves, and establish a related property.

%\bcond\label{PropertyD}
% Suppose that $x_1,$ $x_2$ and $x_3$ are points in $(G,|\cdot|)$, and
% $\alpha_{23}$ denotes a curve in $G$ with endpoints $x_2$ and $x_3$.
%  We say that the quadruple $[x_1,x_2,x_3;\alpha_{23}]$ satisfies {\it Condition \Ref{PropertyD}} if there exists a point
%  $u_1\in\mathbb{S}(x_2, |x_2-x_1|)\cap \alpha_{23}$ such that $$k_{G}(x_1,u_1)\geq 10.$$\medskip
%\econd
%%%%%%%%%%%%%%%%%%%%%%%%%%%%%%%%%%%%%%%%%%%%

%The following result is an immediate consequence of the proof of \cite[Lemma 3.6]{HRWZ} (i.e., Lemma \Ref{H-W-Zh-X-1} as above).

%\begin{lem}\label{sun-2}
% Suppose that the quadruple $[x_1,x_2,x_3;\alpha_{23}]$ satisfies {\it Condition \Ref{PropertyD}}. Then there exists a $100$-quasigeodesic ray $\alpha_{1,1}$ in $G$ starting from $x_1$
%and ending at $x_{1,1}\in \partial_{\sigma}G$ such that for any $x\in\alpha_{1,1}$,
%$$\frac{1}{42}\sigma(x_2,x_1)< \sigma(x_2,x)\leq 5\sigma(x_2,x_1)\;\;\mbox{and}\;\;\ell(\alpha[x,x_{1,1}))\leq 5\delta_{G}(x).$$
%\end{lem}
%%%%%%%%%%%%%%%%%%%%%%%%%%%%%%%%%%%%%%%%%%%%

\bcond\label{PropertyH}
Suppose that $x_1$ and $x_2$ are points in a rectifiably connected, locally compact and non-complete metric space $(X,d)$, $x_3\in \partial_{\sigma}X$, and $\alpha_{23}$ is a ray in $X$ starting from $x_2$ and ending at $x_3$. We say that the quadruple $[x_1,x_2,x_3;\alpha_{23}]$ satisfies {\it Condition \ref{PropertyH}} if
\ben
\item\label{Final-10(1)}
$\sigma(x_1,x_2)\geq M_1^{\frac{1}{13}}\delta_X(x_2)$, and
\item\label{Final-10(2)}
 $\ell(\alpha_{23})\leq 5\delta_X(x_2)$.
\een
\econd
%%%%%%%%%%%%%%%%%%%%%%%%%%%%%%%%%%%%%%%%%%%%

\begin{lem}\label{Wedn-12-7} If the quadruple $[x_1,x_2,x_3;\alpha_{23}]$ satisfies Condition \ref{PropertyH},
then $$k_X(x_1,\alpha_{23})\geq \frac{1}{13}\log M_1-\log 5.$$ \end{lem}
\bpf
By \eqref{eq-neq-1}, we see that for any $x\in \alpha_{23}$,
$$k_X(x_1,x) \geq \log\Big(1+\frac{\sigma(x_1,x)}{\delta_X(x)}\Big).$$

Since $\sigma(x_1,x)\geq \sigma(x_1,x_2)-\ell(\alpha_{23})$ and $\delta_X(x)\leq \ell(\alpha_{23})$, we obtain from the assumptions in Condition \ref{PropertyH} that
$$k_X(x_1,x) \geq \frac{1}{13}\log M_1-\log 5,$$
 as required.
\epf

%%%%%%%%%%%%%%%%%%%%%%%%%%%%%%%%%%%%%%%%%%%
%%%%%%%%%%%%%%%%%%%%%%%%%%%%%%%%%%%%%%%%%%%
\section{Diameter version of the quasiconvexity condition with respect to the inner metric}\label{sec-5}
%%%%%%%%%%%%%%%%%%%%%%%%%%%%%%%%%%%%%%%%%%%
%%%%%%%%%%%%%%%%%%%%%%%%%%%%%%%%%%%%%%%%%%%

The purpose of this section is to prove that, under Condition \ref{PropertyB}, the image of a quasihyperbolic geodesic in $G'$ under $f^{-1}$ satisfies a diameter version of the quasiconvexity condition with respect to the inner metric.

 For convenience, in the following, we write $z'=f(z)$ for points $z$ in $G$ and $\beta'=f(\beta)$ for curves $\beta$ in $G$.

Assume that $G$ satisfies Condition \ref{PropertyB}.
Since $G'$ is $M$-uniform, we know that the identity mapping from $(G',d')$ to $(G',\sigma(d'))$ is a homeomorphism, and thus, \cite[Proposition 2.8]{BHK} implies that $(G', k_{G'})$ is geodesic.
Now, we are going to prove the following theorem.

\begin{thm}\label{sat-11}
Let $x_1$, $x_2$ be two points in $G$, and let $\gamma\subset G$ be a curve joining $x_1$ and $x_2$  such that its image $\gamma'$ is a quasihyperbolic geodesic in $G'$. Then $\diam(\gamma)\leq B_0\sigma(x_1,x_2)$, where $B_0$ is from Theorem \Ref{Thm-HRWZ}.
\end{thm}
\bpf
We prove the theorem by considering two cases.
For the case when there exists a point $\xi_0\in\gamma$ such that
$\delta_{G}(\xi_0)> 2\max\{\sigma(x_1,\xi_0),$ $\sigma(x_2,\xi_0)\}$, we claim that

\be\label{lem-6.1}
\ell(\gamma)\leq 3e^{2M^2}M^2\sigma(x_1,x_2). \ee

 The assumption in this case ensures that the segment $[x_1,x_2]\subset \mathbb{B}(\xi_0, \frac{1}{2}\delta_G(\xi_0))$. Hence we get that for any $x\in[x_1,x_2]$,
\beqq\label{21-01-1}
\delta_{G}(x)\geq \frac{1}{2}\delta_{G}(\xi_0)\geq \frac{1}{2}|x_1-x_2|=\frac{1}{2}\sigma (x_1,x_2),
\eeqq
and thus, we get
\be\label{21-01-2}
k_G(x_1,x_2)\leq \int_{[x_1,x_2]}\frac{|dx|}{\delta_{G}(x)}\leq \frac{2\sigma(x_1,x_2)}{\delta_{G}(\xi_0)}  \leq 2.
\ee

Since $f$ is $M$-quasihyperbolic, we see that its inverse $f^{-1}$ is also $M$-quasihyperbolic. Thus it follows from Lemma \Ref{NY-1} that $\gamma$ is $M^2$-quasigeodesic, which gives that
\be\label{21-01-3}
\ell_{k}(\gamma)\leq  M^2 k_G(x_1,x_2).
\ee
Then \eqref{vai-lem3.2} leads to
\be\label{21-01-4}
\log\Big(1+\frac{\ell(\gamma)}{\delta_{G}(x_1)}\Big)\leq \ell_{k}(\gamma)\leq  2M^2,
\ee
 and since $\log(1+t)\geq e^{-2M^2}t$ for $0<t\leq e^{2M^2}-1$, we have
\be\label{21-01-5}
e^{-2M^2}\frac{\ell(\gamma)}{\delta_G(x_1)}\leq\log\Big(1+\frac{\ell(\gamma)}{\delta_{G}(x_1)}\Big).
 \ee

Now, we conclude from the inequalities \eqref{21-01-2}$-$\eqref{21-01-5} that
$$
e^{-2M^2}\frac{\ell(\gamma)}{\delta_G(x_1)}\leq 2 M^2\frac{\sigma(x_1,x_2)}{\delta_{G}(\xi_0)},
$$
 which, together with the fact $\delta_{G}(x_1)\leq \delta_G(\xi_0)+|x_1-\xi_0|\leq\frac{3}{2}\delta_{G}(\xi_0)$, implies that
\be\label{25-1-13-1}
\ell(\gamma)\leq 3e^{2M^2}M^2\sigma(x_1,x_2).\ee

\begin{figure}
\begin{center}
\includegraphics[width=11cm]{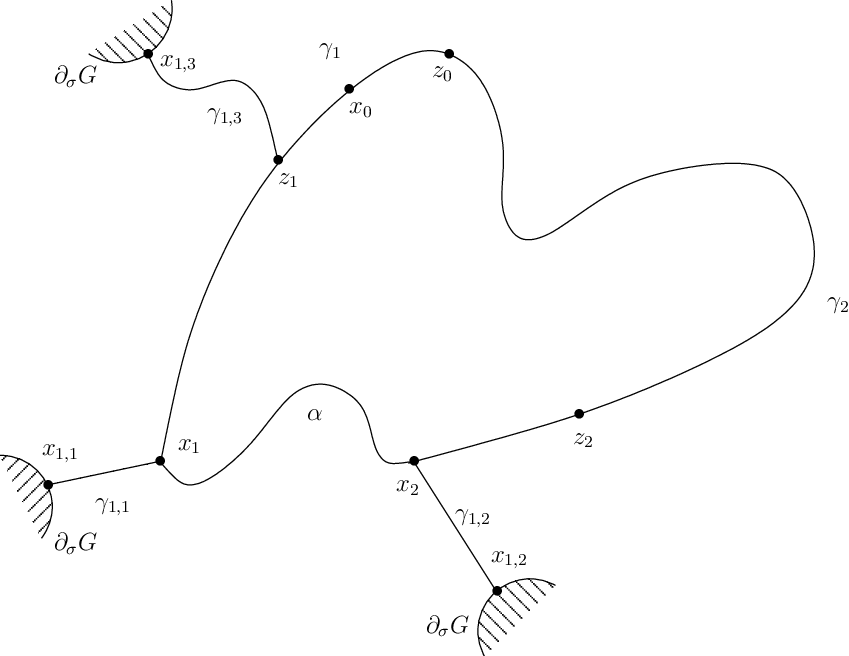}
\caption{The related points and curves.}
\label{fig8}
\end{center}
\end{figure}

For the remaining case, that is, for any $x\in\gamma$,
$\delta_{G}(x)\leq 2\max\{\sigma(x_1,x),\sigma(x_2,x)\},$
let $z_0$ be a point in $\gamma$ such that
\beq\label{chr-2-3}
\sigma(x_1,z_0)=\frac{1}{2}\diam(\gamma)\eeq $($see Figure~\ref{fig8}$)$. Then the following are four lemmas concerning $z_0$.

\blem\label{Corr-12-Tues-2} Suppose that $\diam(\gamma)\geq B_0\sigma(x_1,x_2)$. Then
 $$\min\{\sigma(x_1,z_0), \sigma(x_2,z_0)\}\geq \frac{1}{3}\max\big\{2\sigma(x_1,z_0), B_0\sigma(x_1,x_2)\big\}.$$\elem
\bpf
It follows from \eqref{chr-2-3} and the assumption of the lemma that $$\sigma(x_1,x_2)\leq 2B_0^{-1}\sigma(x_1,z_0).$$
Then, by the triangle inequality $$\sigma(x_2,z_0)\geq \sigma(x_1,z_0)-\sigma(x_1,x_2),$$ we see that the lemma is true.
\epf

\blem\label{changsha-12} Suppose that $\diam(\gamma)\geq B_0\sigma(x_1,x_2)$ and $\max\{\delta_G(x_1),\delta_G(x_2)\}\leq 2\sigma(x_1,x_2)$. Then
 $$\max\{\delta_G(x_1),\delta_G(x_2)\}\leq 6B_0^{-1}\min\{\sigma(x_1,z_0), \sigma(x_2,z_0)\}.$$\elem
\bpf
Obviously, the lemma follows from Lemma \ref{Corr-12-Tues-2}.
\epf

\blem\label{chr-2} Suppose that $(1)$ $\diam(\gamma)\geq B_0\sigma(x_1,x_2)$; $(2)$ $\max\{\delta_G(x_1),\delta_G(x_2)\}\leq 2\sigma(x_1,x_2)$, and
$(3)$ there are $i\in\{1,2\}$ and a point $u\in\gamma$ such that $\sigma(x_i,u)\geq 6B_0^{-1}\sigma(x_1,z_0)$. Then
 $$\delta_G(u)\leq 2\sigma(x_i,u).$$\elem
\bpf
Assume that $i\in\{1,2\}$ and $u\in\gamma$. If $\sigma(x_i,u)\geq 6B_0^{-1}\sigma(x_1,z_0)$, then we infer from
 Lemma \ref{changsha-12} that
$$\delta_G(u)\leq \delta_G(x_i)+\sigma(x_i,u)\leq  6B_0^{-1}\sigma(x_1,z_0)+\sigma(x_i,u)\leq 2\sigma(x_i,u),$$
which is the required estimate.
\epf

\blem\label{chr-1} Suppose that $(1)$ $\diam(\gamma)\geq B_0\sigma(x_1,x_2)$; $(2)$ $\max\{\delta_G(x_1),\delta_G(x_2)\}\leq 2\sigma(x_1,x_2)$, and
$(3)$ there are $i\in\{1,2\}$ and a point $v_i\in \gamma[x_i, z_0]$ such that
$$\min\{ \sigma(x_1,v_i), \sigma(x_2,v_i)\}\leq M_1^{-\frac{1}{4}}\sigma(x_1,z_0).$$  Then for all $v\in\gamma\backslash\gamma[x_i,z_0)$,
  $$k_G(v_i,v)\geq \frac{1}{4}M^{-2}\log M_1-1.$$ \elem
\bpf
We start the proof with the check of the following inequality
\be\label{shantou-3}
\sigma(x_1,v_i)\leq 2M_1^{-\frac{1}{4}}\sigma(x_1,z_0).
\ee

For the case $\min\{ \sigma(x_1,v_i), \sigma(x_2,v_i)\}=\sigma(x_1,v_i)$, the inequality is obvious. For the remaining case, that is,
 $\min\{ \sigma(x_1,v_i), \sigma(x_2,v_i)\}=\sigma(x_2,v_i)$, we obtain from Lemma \ref{Corr-12-Tues-2} that
$$\sigma(x_1,v_i)\leq \sigma(x_2,v_i)+\sigma(x_1,x_2)< 2M_1^{-\frac{1}{4}}\sigma(x_1,z_0),$$
as required.\medskip

Since Lemma \ref{Corr-12-Tues-2} and \eqref{shantou-3} guarantee
$$\sigma(x_1,z_0)-\sigma(x_1,v_i)\geq \big(\frac{2}{3}-2M_1^{-\frac{1}{4}}\big) \sigma(x_1,z_0),$$
and because Lemma \ref{changsha-12} and \eqref{shantou-3} ensure
$$\sigma(x_1,v_i)+\delta_G(x_1)\leq (6B_0^{-1}+2M_1^{-\frac{1}{4}}\big)\sigma(x_1,z_0), $$
we have that
\begin{eqnarray*}
k_G(v_i,v)\nonumber&\geq& M^{-2}\ell_{k}(\gamma[v_i,v]) \;\;\;\;\;\;\;\;\;\;\;\;\;\;\;\;\;\;\;\;\;\;\;\;\;\,\;\;\;\;\;\;\;\;\;\;\;\;\;\;\;\;  (\mbox{by Lemma \Ref{NY-1}})
\\ \nonumber&\geq& M^{-2}k_G(v_i,z_0)\geq M^{-2}\log\Big(1+\frac{\sigma(v_i,z_0)}{\delta_G(v_i)}\Big)\;\;\; (\mbox{by \eqref{eq-neq-1}})
\\ \nonumber&\geq& M^{-2}\log\Big(1+\frac{\sigma(x_1,z_0)-\sigma(x_1,v_i)}{\sigma(x_1,v_i)+\delta_G(x_1)}\Big)
\\ \nonumber&>& \frac{1}{4}M^{-2}\log M_1-1, \;\;\;\;  (\mbox{by Remark \Ref{nw-1} (\ref{HWRZ-25-4})} )\end{eqnarray*}
which completes the proof.
\epf

Now, we are ready to state and prove the following relation.

\bprop\label{lem-6.2}
$\diam(\gamma)< B_0\sigma(x_1,x_2).$
\eprop
\bpf
It follows from the assumption of this case that
\be\label{Wednes-12-22-1} \max\{\delta_G(x_1),\delta_G(x_2)\}\leq 2\sigma(x_1,x_2).\ee

For each $i\in \{1,2\}$,
we choose $x_{1,i}\in \partial_{\sigma} G$ so that
\be\label{Wedns-12-7}
\delta_{G}(x_i)=|x_i-x_{1,i}|, 
\ee
and let $\gamma_{1,i}=[x_i,x_{1,i})=[x_i,x_{1,i}]\setminus \{x_{1,i}\}$,
where $[x_i,x_{1,i}]$ denotes the segment in $\IR^n$ with endpoints
$x_i$ and $x_{1,i}$ (see Figure~\ref{fig8}).

We prove the proposition by contradiction. Suppose, to the contrary, that
\be\label{Exter-1}\diam(\gamma)\geq B_0\sigma(x_1,x_2).\ee
Then we infer from \eqref{Wednes-12-22-1} that all assumptions in Lemmas \ref{Corr-12-Tues-2} and \ref{changsha-12} are satisfied.

Let $\alpha$ be a curve in $G$ connecting $x_1$ and $x_2$ with
$$\ell(\alpha)<\frac{5}{4}\sigma(x_1,x_2),$$ and for $i\in \{1,2\}$, let
\be\label{12-Tues-0} \gamma_i=\alpha\cup\gamma[x_i,z_0].\ee
Then it follows from Lemma \ref{Corr-12-Tues-2} that
\be\label{changsha-11}\ell(\alpha)<\frac{15}{4}B_0^{-1}\min\{\sigma(x_1,z_0), \sigma(x_2,z_0)\}.\ee

Let $x'_0\in \gamma'$ bisect $\gamma'$. Without loss of generality, we assume that $x'_0\in \gamma'[x'_1,z'_0]$.
 We consider two possibilities:\medskip

\noindent {\bf Subcase 5.2.1}\label{Zhou-H-1} {\it
Assume that $\sigma(x_1,x_0)\geq M_1^{-1}\sigma(x_1,z_0).$  }
\medskip

 It follows from \eqref{changsha-11} that (see Remark \ref{nw-1} (\ref{HWRZ-25-4}))$$\sigma(x_1,x_2)\leq \ell(\alpha)<\frac{15B_0^{-1}}{4}\sigma(x_1,z_0)<M_1^{-1}\sigma(x_1,z_0).$$
Then the assumption of this subcase implies that there exist points $z_1\in \gamma[x_1,x_0]$ and $z_2\in \gamma[x_2,z_0]$ such that
\be\label{Corr-12-Tues-3}\sigma(x_1,z_1)=\sigma(x_1,z_2)=M_1^{-2}\sigma(x_1,z_0)\ee  (see Figure~\ref{fig8}).
Then we deduce from Lemma \ref{changsha-12} that for each $i\in \{1,2\}$,
\be\label{Wedns-12-8}\sigma(x_1,z_i)=M_1^{-2}\sigma(x_1,z_0)\geq \frac{1}{6}B_0M_1^{-2}\delta_G(x_i).\ee

The following two lemmas concern the point $z_1$ and the image $z'_2$.
\blem\label{Zhou-W-0} For any $z\in \gamma_2$,
$ k_G(z_1,z)>\frac{1}{4}M^{-2}\log M_1-1.$
\elem
\bpf
Since $\gamma_2=\alpha\cup\gamma[x_2,z_0]$, we may divide the proof of the lemma into two cases. For the first case when $z\in\alpha$, it follows from \eqref{changsha-11} and \eqref{Corr-12-Tues-3} that
$$\sigma(x_1,z_1)-\ell(\alpha)\geq \big(M_1^{-2}-\frac{15}{4}B_0^{-1}\big)\sigma(x_1,z_0),$$
and from Lemma \ref{changsha-12} and \eqref{changsha-11} that
$$\delta_G(x_1)+\ell(\alpha)\leq \frac{39}{4}B_0^{-1}\sigma(x_1,z_0).$$
Thus we get
$$ k_G(z_1,z) \geq \log\Big(1+\frac{\sigma(z_1,z)}{\delta_G(z)}\Big)>\log\Big(1+\frac{\sigma(x_1,z_1)-\ell(\alpha)}{\delta_G(x_1)+\ell(\alpha)}\Big)\\ >\frac{3}{2}\log M_0.$$

For the remaining case, that is, $z\in\gamma[x_2,z_0]$, we note that (\ref{Corr-12-Tues-3}) ensures the following
$$\min\{\sigma(x_1,z_1),\sigma(x_2,z_1)\}\leq M_1^{-2}\sigma(x_1,z_0)$$
by letting $i=1$, $v_1=z_1$ and $v=z$ in Lemma \ref{chr-1}. Thus we see that
the required inequality follows from Lemma \ref{chr-1}.
\epf

\blem \label{chr-3}
$\delta_{G'}(z'_2)\leq M_2^{-1}\delta_{G'}(x'_0).$
\elem
\bpf
Since (\ref{Corr-12-Tues-3}) guarantees that Lemma \ref{chr-1} works for the situation where $i=2$, $v_i=z_2$ and $v=x_0$, it follows that
$$k_G(z_2,x_0)\geq \frac{1}{4}M^{-2}\log M_1-1.$$

On the other hand, Lemma \Ref{lab-40} and the facts that $x'_0$ bisects $\gamma'$ and $z'_2\in \gamma'[x'_2,x'_0]$ ensure the following
$$d'(z'_2,x'_0) \leq \nu_2\delta_{G'}(x'_0).$$
Hence we get
\begin{eqnarray*}\frac{1}{4}M^{-2}\log M_1-1 &\leq& k_G(z_2,x_0)\;\;\;\; \\ \nonumber&\leq&  M k_{G'}(z'_2,x'_0)\;\;\;\; (\mbox{since $f|_G$ is $M$-quasihyperbolic})
\\ \nonumber&\leq& 4M^3\log\Big(1+\frac{\nu_2\delta_{G'}(x'_0)}{\min\{\delta_{G'}(z'_2),\delta_{G'}(x'_0)\}}\Big),\;\;\;\,\;\;\;\;\;\;(\mbox{by Lemma \Ref{h-W-1}\eqref{h-W-1(1)}})
\end{eqnarray*}
and therefore, $$\delta_{G'}(z'_2)< M_2^{-1}\delta_{G'}(x'_0).$$ This is what we need.
\epf

Next, we apply the assumption that $f$ is $\eta$-quasisymmetric on $\partial_{\sigma}G$ of Condition \ref{PropertyB} to continue the discussions. Before the application, we need to find a point from $\partial_{\sigma}G$, which is determined in the following lemma.

\blem\label{Corec-1}
There exist a point $x_{1,3}\in \partial_{\sigma} G$ and a $100$-quasigeodesic ray $\gamma_{1,3}$ in $G$ starting from $z_1$ and ending at $x_{1,3}$
such that for any $x\in\gamma_{1,3}$,
\beqq\frac{1}{42}\sigma(x_1,z_1)\leq\sigma(x_1,x)\leq 5\sigma(x_1,z_1)\;\;\mbox{and}\;\;\ell(\gamma_{1,3}[x,x_{1,3}))\leq 5\delta_{G}(x)\eeqq
  $($see Figure~\ref{fig8}$)$.
\elem
\bpf
It follows from Lemma \Ref{H-W-Zh-X-1} that we only need to check that the quadruple $[z_1,x_1,z_0;\gamma_2]$ satisfies Condition A.
This immediately follows from (\ref{12-Tues-0}), (\ref{Corr-12-Tues-3}) and Lemma \ref{Zhou-W-0}, and so the lemma is proved.
\epf

To reach a contradiction,  an upper bound for the ratio $\sigma(x_{1,1},x_{1,2})/\sigma(x_{1,1},x_{1,3})$ and
a lower bound for the ratio $d'(x'_{1,1},x'_{1,2})/d'(x'_{1,1},x'_{1,3})$ are still required. First, we establish the following upper bound for the ratio $\sigma(x_{1,1},x_{1,2})/\sigma(x_{1,1},x_{1,3})$.
\blem\label{wed-14}
$\sigma(x_{1,1},x_{1,2})\leq 5M_0^{-\frac{3}{2}}\sigma(x_{1,1},x_{1,3})$.
\elem
\bpf
Since
\beqq\label{Zhou-W-H-2}
\sigma(x_{1,1},x_{1,2})\leq \sigma(x_1,x_{1,1})+\sigma(x_1,x_2)+\sigma(x_2,x_{1,2})\leq 5\sigma(x_1,x_2) \;\;\;(\mbox{by (\ref{Wednes-12-22-1})})
\eeqq
and by Remark \ref{nw-1} (\ref{HWRZ-25-4}), $B_0>128M_0^\frac{3}{2} M_1^6+128M_1^2$, we have
\beqq\label{Zhou-W-H-3}
\sigma(x_{1,1},x_{1,3})&\geq & \sigma(x_1,x_{1,3})-\sigma(x_1,x_{1,1}) 
\\ \nonumber&\geq & \frac{1}{42}\sigma(x_1,z_1)-2\sigma(x_1,x_2)\;\;\;\;\;\;\;(\mbox{by (\ref{Wednes-12-22-1}), (\ref{Wedns-12-7}) and Lemma \ref{Corec-1}}) \\ \nonumber
&>& \frac{M_1^{-2}}{42}\sigma(x_1,z_0)-2\sigma(x_1,x_2), \;\;\;\;(\mbox{by (\ref{Corr-12-Tues-3})})\\ \nonumber
&>& M_0^{\frac{3}{2}}\sigma(x_1,x_2),  \;\;\;\;\;\;\;\;\;\;\;\;\;\;\;\;\;\;\;\;\;\;\;\;(\mbox{by Lemma \ref{Corr-12-Tues-2}})
\eeqq
we see that the lemma holds.
\epf

To get a lower bound for the ratio $d'(x'_{1,1},x'_{1,2})/d'(x'_{1,1},x'_{1,3})$, we need the following three chains of inequalities.

\blem \label{wed-16} \ben
\item\label{wed-16(1)}
$\delta_{G'}(z'_1)\leq d'(z'_1,x'_{1,3})\leq \ell(\gamma'_{1,3})\leq M_2^{-1}\delta_{G'}(x'_0).$
\item\label{wed-16(2)}
$\delta_{G'}(x'_1)\leq d'(x'_1,x'_{1,1})\leq \ell(\gamma'_{1,1}) \leq M_2^{-1}\delta_{G'}(z'_1)$.
\item\label{wed-16(3)}
$\delta_{G'}(x'_2)\leq d'(x'_2,x'_{1,2})\leq \ell(\gamma'_{1,2})\leq M_2^{-1}\delta_{G'}(z'_2).$
\een
\elem
\bpf
We apply Lemmas \Ref{H-W-Zhou-1} and \ref{Wedn-12-7} to prove this lemma.
First, since $M_1^{-2}\geq 6B_0^{-1}$, we know from (\ref{Wednes-12-22-1}), (\ref{Exter-1}) and (\ref{Corr-12-Tues-3}) that Lemma \ref{chr-2} works for the situation where $i=1$ and $u=z_1$. It follows that
$$\delta_G(z_1)\leq 2\sigma(x_1,z_1),$$
and thus, (\ref{Corr-12-Tues-3}) and the assumption of the subcase guarantee that
\be\label{Zhou-W-1}
\sigma(x_0,z_1)\geq \sigma(x_1,x_0)-\sigma(x_1,z_1) \geq (M_1-1)\sigma(x_1,z_1)\geq \frac{1}{2}(M_1-1)\delta_G(z_1).\ee
Now, Lemma \ref{Corec-1} ensures that the quadruple $[x_0,z_1,x_{1,3};\gamma_{1,3}]$ satisfies Condition \ref{PropertyH}. By Lemma \ref{Wedn-12-7}, we have
$$k_{G'}(x'_0,\gamma'_{1,3})\geq M^{-1}k_{G}(x_0,\gamma_{1,3})> \big(\frac{1}{13}\log M_1-\log 5\big)M^{-1},$$
which, together with Lemma \Ref{lab-40} and the inequality
$$d'(z'_1,x'_0)\leq \nu_2\delta_{G'}(x'_0)$$
(since $G'$ is uniform and $x'_0$ bisects $\gamma'$),
shows that the quadruple $[x'_0, z'_1, x'_{1,3};\gamma'_{1,3}]$ satisfies Condition B.
Then we conclude from Lemma \Ref{H-W-Zhou-1} that
\beqq\label{Zhou-W-3} d'(z'_1,x'_{1,3})\leq \ell(\gamma'_{1,3})\leq M_2^{-1}\delta_{G'}(x'_0),\eeqq
which demonstrates that the statement \eqref{wed-16(1)} of the lemma is true.

Before the continuation of the proof, we prove the following estimates for the quantities $\sigma(x_1,z_0)$ and $\sigma(x_2,z_2)$ in terms of $\delta_G(x_2)$:
\beq\label{Wedns-12-10} 2M_1^{-2}\sigma(x_1,z_0)&>&\sigma(x_1,z_2)+\sigma(x_1,x_2)\,\;\;\;(\mbox{by Lemma \ref{Corr-12-Tues-2} and (\ref{Corr-12-Tues-3})})
\\ \nonumber&\geq& \sigma(x_2,z_2)\\ \nonumber&\geq&  \sigma(x_1,z_2)-\sigma(x_1,x_2)
\\ \nonumber&\geq& (1-\frac{3 M_1^2}{B_0}) \sigma(x_1,z_2)\;\;\;\;\,(\mbox{by Lemma \ref{Corr-12-Tues-2} and \eqref{Corr-12-Tues-3}})
\\ \nonumber&>& M_0^{\frac{3}{2}}M_1\delta_G(x_2).\;\;\;\;\;\;\;\;\;\;\;\;\,(\mbox{by \eqref{Wedns-12-8}})\eeq

Now, we are ready to complete the proof.
By replacing (\ref{Zhou-W-1}), Lemma \ref{Corec-1}, $[x_0,z_1,x_{1,3};\gamma_{1,3}]$, $\gamma'_{1,3}$ and $[x'_0, z'_1, x'_{1,3};\gamma'_{1,3}]$ with
(\ref{Wedns-12-8}), (\ref{Wedns-12-7}), $[z_1,x_1,x_{1,1};\gamma_{1,1}]$, $\gamma'_{1,1}$ and $[z'_1, x'_1, x'_{1,1}; \gamma'_{1,1}]$ (resp. (\ref{Wedns-12-10}), (\ref{Wedns-12-7}), $[z_2, x_2, x_{1,2}; \gamma_{1,2}]$, $\gamma'_{1,2}$ and $[z'_2, x'_2, x'_{1,2};\gamma'_{1,2}]$), respectively, similar arguments as in the
proof of the statement \eqref{wed-16(1)} show that the statement \eqref{wed-16(2)} (resp.  the statement \eqref{wed-16(3)}) of the lemma is true as well.
\epf

The following is our required lower bound for the ratio $d'(x'_{1,1},x'_{1,2})/d'(x'_{1,1},x'_{1,3})$.

\blem\label{wed-15}
$d'(x'_{1,1},x'_{1,2})\geq \frac{1}{4}\nu_2^{-2} M_2 d'(x'_{1,1},x'_{1,3}).$
\elem
\bpf
Since $x'_0$ bisects $\gamma'$, and since $z'_1$ is between $x'_1$ and $x'_0$, we have
$$\min\{\ell(\gamma'[x'_1,z'_1]),\ell(\gamma'[x'_2,z'_1])\}=\ell(\gamma'[x'_1,z'_1]).$$
Then Lemma \Ref{lab-40} gives
\be\label{21-01-22-1}
\ell(\gamma')\leq \nu_2d'(x'_1,x'_2)\;\;\mbox{and}\;\;d'(x'_1,z'_1)\leq \ell(\gamma'[x'_1,z'_1])\leq\nu_2\delta_{G'}(z'_1).
\ee

Moreover, the last two chains of inequalities in Lemma \ref{wed-16} imply
$$d'(x'_1,x'_{1,1})+d'(x'_2,x'_{1,2})\leq M_2^{-1}(\delta_{G'}(z'_1)+\delta_{G'}(z'_2)),$$
and thus, we obtain from Lemmas \ref{chr-3} and \ref{wed-16}\eqref{wed-16(1)} that
\be\label{21-01-22-2}
d'(x'_1,x'_{1,1})+d'(x'_2,x'_{1,2})\leq 2 M_2^{-2}\delta_{G'}(x'_0).
\ee

Secondly, the first two chains of inequalities in Lemma \ref{wed-16} imply
\be\label{21-01-22-3}
\ell(\gamma')=2\ell(\gamma'[x'_1,x'_0])\geq 2(\delta_{G'}(x'_0)-\delta_{G'}(x'_1))\geq 2(1-M_2^{-2})\delta_{G'}(x'_0)
\ee
and
\be\label{21-01-22-4}
d'(x'_1,x'_{1,1})+d'(z'_1,x'_{1,3})\leq M_2^{-1}(\delta_{G'}(z'_1)+\delta_{G'}(x'_0)).
\ee
Hence we conclude that
\beqq\label{Zhou-WW-5}
d'(x'_{1,1},x'_{1,2})&\geq& d'(x'_1,x'_2)-d'(x'_1,x'_{1,1})-d'(x'_2,x'_{1,2})
\\ \nonumber&>& \nu_2^{-1}\ell(\gamma')-2 M_2^{-2}\delta_{G'}(x'_0) \;\;\; \;\;\;\;\;\; \;\;\;\; (\mbox{by \eqref{21-01-22-1} and \eqref{21-01-22-2}})\\ \nonumber
&\geq& (2\nu_2^{-1}(1-M_2^{-2})-2M_2^{-2}) \delta_{G'}(x'_0) \;\;\;\;(\mbox{by \eqref{21-01-22-3}}) \\ \nonumber&>&\frac{1}{2}\nu_2^{-1}\delta_{G'}(x'_0) 
\eeqq
and
\beqq\label{Zhou-W-6}
d'(x'_{1,1},x'_{1,3})&\leq& d'(x'_1,x'_{1,1})+d'(z'_1,x'_{1,3})+d'(x'_1,z'_1)
\\ \nonumber&\leq& M_2^{-1}\big(\delta_{G'}(z'_1)+\delta_{G'}(x'_0)\big)+ \nu_2\delta_{G'}(z'_1)\;\;\;(\mbox{by \eqref{21-01-22-1} and \eqref{21-01-22-4}})
\\ \nonumber&<& 2\nu_2 M_2^{-1}\delta_{G'}(x'_0),\;\;\;\;\;\;\;\;\;\;\;\;\;\;\;\;\;\;\;\;\;\;\;\;\;\;\;\;\;\;\;\;\;(\mbox{by Lemma \ref{wed-16}\eqref{wed-16(1)}})
\eeqq
from which Lemma \ref{wed-15} follows.\epf

Now, we are ready to reach a contradiction. By the assumption \eqref{PropertyB-3} of Condition \ref{PropertyB}, we know that
$$\frac{d'(x'_{1,1},x'_{1,2})}{d'(x'_{1,1},x'_{1,3})}\leq \eta\Big(\frac{\sigma(x_{1,1},x_{1,2})}{\sigma(x_{1,1},x_{1,3})}\Big).$$ Then it follows from Lemmas \ref{wed-14} and \ref{wed-15} that
$$M_2\leq 4\nu_2^2\eta(5M_0^{-\frac{3}{2}}).$$ However, $M_2\geq 10\nu_3^4\eta(\nu_3)$. This obvious contradiction implies
that the inequality in (\ref{Exter-1}) is not true, proving the proposition in this subcase.
\medskip

\noindent{\bf Subcase 5.2.2}\label{Zhou-H-2} {\it
Assume that $\sigma(x_1,x_0)< M_1^{-1}\sigma(x_1,z_0).$ }\medskip

As in the previous subcase, we are going to reach a contradiction by applying the assumption that $f$ is $\eta$-quasisymmetric on $\partial_{\sigma}G$  in Condition \ref{PropertyB}. For this, we need to choose two points  from $\partial_{\sigma}G$ (i.e., the point $u_{1,1}$ in Lemma \ref{Zhou-W-9} and the other one $u_{1,2}$ in Lemma \ref{Zhou-W-8} below).

\begin{figure}
\begin{center}
\includegraphics[width=11cm]{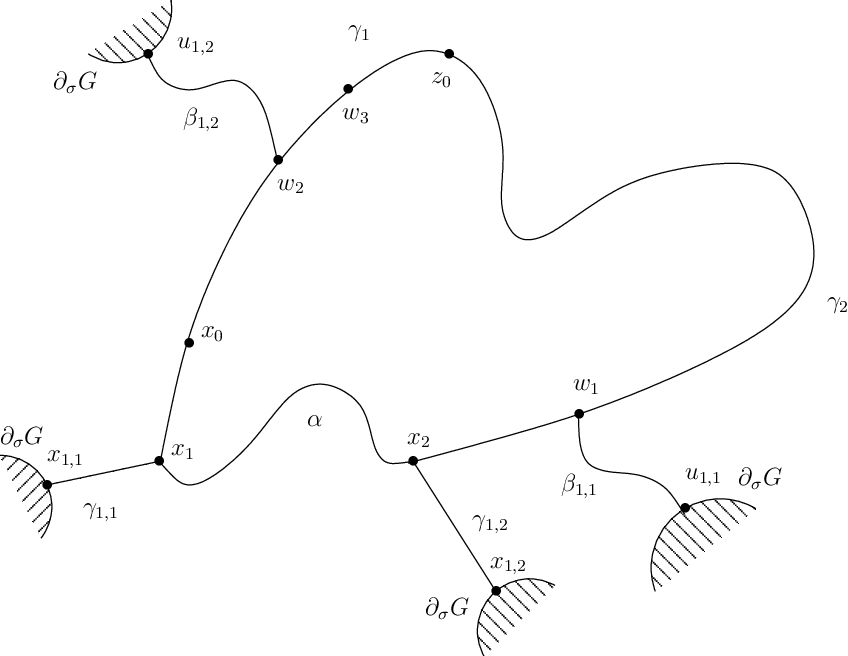}
\caption{The related points and curves.}
\label{fig9}
\end{center}
\end{figure}

First, we choose three points $w_1$, $w_2$ and $w_3$ from $\gamma$ as follows. Since Lemma \ref{Corr-12-Tues-2} shows $$\sigma(x_2,z_0)\geq \frac{2}{3}\sigma(x_1,z_0),$$
it ensures that there exists a point $w_1\in \gamma[x_2,z_0]$ satisfying
\be\label{mon-4} \sigma(x_2,w_1)=2M_1^{-\frac{1}{3}}\sigma(x_1,z_0)\ee  (see Figure~\ref{fig9}). Then by replacing \eqref{Corr-12-Tues-3} with \eqref{mon-4}, a similar reasoning as in the proof of Lemma \ref{Zhou-W-0} shows that the following lemma holds:
\blem\label{changsha-13}
For any $z\in \gamma_1$, $k_G(w_1,z)>\frac{1}{4}M^{-2}\log M_1-1.$
\elem

Moreover, it follows from \eqref{changsha-11} and the assumption of this subcase that
$$\sigma(x_2,x_0)\leq \sigma(x_1,x_0)+\ell(\alpha)<2M_1^{-1}\sigma(x_1,z_0)$$
and
$$\sigma(x_2,z_0)\geq \sigma(x_1,z_0)-\sigma(x_1,x_2)\geq (1-\frac{15}{4}B_0^{-1})\sigma(x_1,z_0),$$
which, together with Remark \ref{nw-1} (\ref{HWRZ-25-4}), shows that 
$$\sigma(x_2,z_0)>\frac{3}{4}\sigma(x_1,z_0)\;\mbox{and}\; 2M_1^{-\frac{5}{12}}\sigma(x_1,z_0)>2M_1^{-\frac{1}{2}}\sigma(x_1,z_0)>M_1^{\frac{1}{2}}\sigma(x_2,x_0).$$
This ensures that there exist $w_2\in\gamma[x_0,z_0]$ and $w_3\in \gamma[w_2,z_0]$ such that
\be\label{Wend-12-12}\sigma(x_2,w_2)=2M_1^{-\frac{1}{2}}\sigma(x_1,z_0)\;\;\mbox{and}\;\;\sigma(x_2,w_3)=2M_1^{-\frac{5}{12}}\sigma(x_1,z_0)\ee
 (see Figure~\ref{fig9}).

Now, we are ready to find the required points from $\partial_{\sigma}G$ by applying Lemma \Ref{H-W-Zh-X-1}.

\blem\label{Zhou-W-9}
There are a point $u_{1,1}\in \partial_{\sigma} G$ and a $100$-quasigeodesic ray $\beta_{1,1}$ in $G$ starting from $w_1$ and ending at $u_{1,1}$ such that for any $x\in\beta_{1,1}$,
$$ \frac{1}{42}\sigma(x_2,w_1)\leq\sigma(x_2,x)\leq 5\sigma(x_2,w_1)\;\;\mbox{and}\;\;\ell(\beta_{1,1}[x,u_{1,1}))\leq 5\delta_{G}(x)$$  $($see Figure~\ref{fig9}$)$.
\elem
\bpf
We conclude from Lemma \ref{Corr-12-Tues-2} and (\ref{mon-4}) that
$$\sigma(x_2,z_0)\geq \frac{1}{3}M_1^{\frac{1}{3}}\sigma(x_2,w_1),$$
which, together with Lemma \ref{changsha-13}, guarantees that the quadruple $[x_2,w_1,z_0;\gamma_1]$ satisfies Condition A. Then our lemma follows from Lemma \Ref{H-W-Zh-X-1}.
\epf

\blem\label{Zhou-W-8}
There exist a point $u_{1,2}\in \partial_{\sigma} G$ and a $100$-quasigeodesic ray $\beta_{1,2}$ in $G$ starting from $w_2$ and ending at $u_{1,2}$
such that for any $x\in\beta_{1,2}$,
$$\frac{1}{42}\sigma(x_2,w_2)\leq\sigma(x_2,x)\leq 5\sigma(x_2,w_2)\;\;\mbox{and}\;\;\ell(\beta_{1,2}[x,u_{1,2}))\leq 5\delta_{G}(x)
$$  $($see Figure~\ref{fig9}$)$.\elem
\bpf
It follows from (\ref{Wend-12-12}) that
$$\min\{ \sigma(x_1,w_2), \sigma(x_2,w_2)\}\leq 2M_1^{-\frac{1}{2}}\sigma(x_1,z_0).$$
Hence Lemma \ref{chr-1} implies $$k_G(w_2,\gamma[x_2,z_0])\geq \frac{1}{4}M^{-2}\log M_1-1.$$
Moreover, we derive from Lemma \ref{Corr-12-Tues-2} and (\ref{Wend-12-12}) that
$$\sigma(x_2,z_0)\geq \frac{1}{3}M_1^{\frac{1}{2}}\sigma(x_2,w_2).$$
These show that the quadruple $[w_2,x_2,z_0;\gamma[x_2,z_0]]$ satisfies Condition A, and so, the lemma follows from Lemma \Ref{H-W-Zh-X-1}.%\medskip
\epf

The following is an upper bound of the ratio $\sigma(u_{1,2},x_{1,2})/\sigma(x_{1,2},u_{1,1})$.

\blem\label{Zhou-W-H-17}
$\sigma(u_{1,2},x_{1,2})\leq 270 M_1^{-\frac{1}{6}}\sigma(x_{1,2},u_{1,1}).$
\elem
\bpf
Since
\beqq
\sigma(x_{1,2},u_{1,1})&\geq& \sigma(x_2,u_{1,1})-\sigma(x_2,x_{1,2})
\\ \nonumber&\geq& \frac{1}{42}\sigma(x_2,w_1)-2\sigma(x_1,x_2)\;\;\;(\mbox{by (\ref{Wednes-12-22-1}) and Lemma \ref{Zhou-W-9}})
\\ \nonumber&>& \frac{1}{45}\sigma(x_2,w_1),\;\;\;\;\;\;\;\;\;\;\;\;\;\;\;\;\;\;\;\;\;(\mbox{by Lemma \ref{Corr-12-Tues-2} and (\ref{mon-4})})\eeqq
and because
\beqq \sigma(u_{1,2},x_{1,2}) \nonumber&\leq&\sigma(x_2,u_{1,2})+\sigma(x_2,x_{1,2})
\\ \nonumber&<& 5\sigma(x_2,w_2)+2\sigma(x_1,x_2)\;\;\;(\mbox{by (\ref{Wednes-12-22-1}) and Lemma \ref{Zhou-W-8}})
\\ \nonumber&<& 6M_1^{-\frac{1}{6}}\sigma(x_2,w_1),\;\;\;\;\;\;\;\;\;\;\;\;\;\;(\mbox{by Lemma \ref{Corr-12-Tues-2}, \eqref{mon-4} and (\ref{Wend-12-12})})
\eeqq we see that the lemma is true.
\epf

We still need a lower bound for the ratio $d'(u'_{1,2},x'_{1,2})/ d'(x'_{1,2},u'_{1,1})$.
Before the statement and the proof of the needed bound, we establish the following two auxiliary lemmas.
\blem\label{wed-17}
\ben
%\item\label{wed-17(1)}
%$\delta_{G'}(w'_2)\leq d'(w'_2,u'_{1,2})\leq \ell{(\beta'_{1,2})} \leq M_2^{-1}\delta_{G'}(x'_0).$
\item\label{wed-17(2)}
$\delta_{G'}(w'_1)\leq d'(w'_1,u'_{1,1})\leq \ell{(\beta'_{1,1})} \leq M_2^{-1}\delta_{G'}(w'_3).$
\item\label{wed-17(3)}
$\delta_{G'}(x'_2)\leq d'(x'_2,x'_{1,2})\leq \ell{(\gamma'_{1,2})} \leq M_2^{-1}\delta_{G'}(w'_1).$
\een
\elem
\bpf We use Lemma \Ref{H-W-Zhou-1} to prove this lemma.
Because (\ref{mon-4}), \eqref{Wend-12-12} and Lemma \ref{Zhou-W-9} ensure that for $z\in\beta_{1,1}$,
$$\sigma(x_2,z)\geq \frac{1}{42}\sigma(x_2,w_1)=\frac{1}{42}M_1^{\frac{1}{12}}\sigma(x_2,w_3),$$
since $B_0\geq 20e^{2\nu_3M_1}$ (see Remark \ref{nw-1}),
and the  combination of (\ref{Wend-12-12}) and Lemma \ref{chr-2} leads to $$\delta_{G}(w_3)\leq 2\sigma(x_2,w_3),$$
we now conclude that for any $z\in\beta_{1,1}$,
\beqq\label{Zhou-W-10}
\sigma(w_3,z)\geq \sigma(x_2,z)-\sigma(x_2,w_3)> \frac{1}{96}M_1^{\frac{1}{12}}\delta_{G}(w_3).
\eeqq
Thus the assumption that $f|_G$ is $M$-quasihyperbolic of Condition \ref{PropertyB} implies that
$$k_{G'}(w'_3, \beta'_{1,1})\geq M^{-1}k_{G}(w_3, \beta_{1,1})\geq M^{-1}\log\Big(1+\frac{\sigma(w_3,\beta_{1,1})}{\delta_G(w_3)}\Big)>\frac{1}{20}\nu_2^{-1}\log M_1,$$
which, in conjunction with the fact that $\beta'_{1,1}$ is $\nu_2$-uniform (by Lemma \Ref{shantou-5}) and the inequality
$$d'(w'_1, w'_3)\leq C\delta_{G'}(w'_3)$$
(by Lemma \Ref{lab-40} and the facts that $x_0'$ bisects $\gamma'$ and $w_1'\in \gamma'[x_2',w_3']\subset \gamma'[x_2',x_0']$), shows that
the quadruple $[w'_3, w'_1, u'_{1,1}; \beta'_{1,1}]$ satisfies Condition B.
By Lemma \Ref{H-W-Zhou-1}, we see that the first statement of the lemma holds.

To prove the second chain of inequalities, recall that $x_{1,2}$ is a point on $\partial_{\sigma} G$ such that $\delta_{G}(x_2)=|x_2-x_{1,2}|$ and
$\gamma_{1,2}=[x_2,x_{1,2})$ (see Figure~\ref{fig9}).
Then for any $z\in\gamma_{1,2}$,  we get
\begin{eqnarray} k_G(w_1,z)\nonumber&\geq& \log\Big(1+\frac{\sigma(w_1,z)}{\delta_G(z)}\Big)>\log\Big(1+\frac{\sigma(x_2,w_1)-\delta_G(x_2)}{\delta_G(x_2)}\Big)\\
\nonumber&>&\frac{3}{2}\log M_0,\;\;\;\;\;\;\;(\mbox{by Lemma \ref{changsha-12} and (\ref{mon-4})})\end{eqnarray}
and therefore, $$k_{G'}(w'_1,\gamma'_{1,2})\geq M^{-1}k_G(w_1,\gamma_{1,2})>\frac{3}{2}M^{-1}\log M_0.$$
Thus we see from the fact that $\gamma'_{1,2}$ is $\nu_2$-uniform (by Lemma \Ref{shantou-5}) and the inequality $$d'(x'_2, w'_1)\leq \nu_2\delta_{G'}(w'_1)$$
 (by Lemma \Ref{lab-40} and the facts that $x_0'$ bisects $\gamma'$ and that $w_1'\in \gamma'[x_2',x_0']$) that
the quadruple $[w'_1, x'_2, x'_{1,2}; \gamma'_{1,2}]$ satisfies Condition B.
Then the second chain of inequalities of the lemma follows from Lemma \Ref{H-W-Zhou-1} as well.
\epf

\blem\label{Hz-W-Zhou-1} $d'(u'_{1,2}, u'_{1,1})>4\nu_2d'(u'_{1,1},w'_1).$\elem
\bpf
First, we need a lower bound for the quantity $k_{G}(w_3,\beta_{1,2})$. We derive this bound by applying Lemma \ref{Wedn-12-7}. For this, we verify that the quadruple $[w_3, w_2, u_{1,2}; \beta_{1,2}]$ satisfies Condition \ref{PropertyH}.
Since
\beqq
\sigma(w_3,w_2)&\geq& \sigma(x_2, w_3)-\sigma(x_2,w_2)
\geq (M_1^{\frac{1}{12}}-1)\sigma(x_2,w_2) \;\;\;(\mbox{by (\ref{Wend-12-12})})\\ \nonumber
&>& \frac{1}{2}(M_1^{\frac{1}{12}}-1) \delta_{G}(w_2), \;\;\;\;\;\;\;\;\;\;\;\;\;\;\;\;\;\;\;\;\;\;\;\;\;\;\;\;\;\;\;\;\;\;\;\;\;\,\;\;(\mbox{by  Lemma \ref{chr-2}})
\eeqq
we see from Lemma \ref{Zhou-W-8} that the quadruple $[w_3, w_2, u_{1,2}; \beta_{1,2}]$ satisfies Condition \ref{PropertyH}.
Then Lemma \ref{Wedn-12-7} shows that
\be\label{Zhou-W-14} k_{G}(w_3,\beta_{1,2})> \frac{1}{13}\log M_1-\log 5.\ee

Now, we are ready to prove the lemma. Assume, to the contrary, that
\be\label{shantou-2} d'(u'_{1,2}, u'_{1,1})\leq 4\nu_2d'(u'_{1,1},w'_1).\ee

Since it follows from Lemma \ref{wed-17}\eqref{wed-17(2)} that
$$d'(u'_{1,1},w'_1)\leq M_2^{-1}\delta_{G'}(w'_3)$$
and
$$d'(w'_3,w'_1)\geq \delta_{G'}(w'_3)-\delta_{G'}(w'_1)\geq (1-M_2^{-1})\delta_{G'}(w'_3),$$
we get
$$d'(u'_{1,1},w'_1)\leq M_2^{-1}\delta_{G'}(w'_3)\leq \frac{1}{M_2-1}d'(w'_3,w'_1),$$
and so, \eqref{shantou-2} implies
\be\label{changsha-14}
d'(w'_1,u'_{1,2})\leq d'(u'_{1,1},w'_1)+d'(u'_{1,2}, u'_{1,1})< 5\nu_2M_2^{-1}d'(w'_3,w'_1).
\ee
Then we  conclude from Lemma \Ref{lab-40} that
\be\label{changsha-15}
d'(w'_1,u'_{1,2})\leq 5\nu_2M_2^{-1}\ell(\gamma'[w'_1,w'_2])< 5\nu_2^2M_2^{-1}d'(w'_1,w'_2),
\ee
and thus,
\be\label{changsha-16}
d'(w'_2,u'_{1,2})\geq d'(w'_1,w'_2)-d'(w'_1,u'_{1,2})\geq \big(\frac{1}{5}\nu_2^{-2}M_2-1\big )d'(w'_1,u'_{1,2}).
\ee
Since $\gamma'$ is a quasihyperbolic geodesic, $\beta'_{1,2}$ is $\nu_2$-uniform (by Lemma \Ref{shantou-5}) and $w_2'\in \gamma'[w_1',x_0']$, we see from Lemma \Ref{lab-40} and \eqref{changsha-14}$-$\eqref{changsha-16} that
the sextuple $[w'_1,w'_2,w'_3,$ $u'_{1,2};\gamma'[w'_1,w'_2],\beta'_{1,2}]$
satisfies Condition C. Then we derive from Lemma \Ref{H-W-2-zhou} that
$$k_{G}(w_3,\beta_{1,2})\leq Mk_{G'}(w'_3,\beta'_{1,2})\leq 12M^2\log \nu_3,$$
which contradicts (\ref{Zhou-W-14}). Hence the lemma is proved.
\epf

Our required lower bound for the ratio $d'(u'_{1,2},x'_{1,2})/d'(x'_{1,2},u'_{1,1})$ is as follows.

\blem\label{imp-1}
$d'(u'_{1,2},x'_{1,2})>\frac{1}{2}d'(x'_{1,2},u'_{1,1})$.
\elem
\bpf
Since Lemma \ref{wed-17}\eqref{wed-17(3)} gives $$d'(x'_2,w'_1)\geq \delta_{G'}(w'_1)-\delta_{G'}(x'_2)\geq (1-M_2^{-1})\delta_{G'}(w'_1),$$
again, we obtain from Lemma \ref{wed-17}\eqref{wed-17(3)} that
$$d'(x'_2,x'_{1,2})\leq M_2^{-1}\delta_{G'}(w'_1)\leq \frac{1}{M_2-1}d'(x'_2,w'_1).$$
Hence we have that
\beq\label{Add25-1} d'(x'_{1,2},u'_{1,1})&\leq&  d'(x'_2,x'_{1,2})+d'(x'_2,w'_1)+d'(u'_{1,1},w'_1)
\\ \nonumber&<&(\frac{1}{M_2-1}+1)d'(x'_2,w'_1)+d'(u'_{1,1},w'_1).\eeq
Morever, since $x'_0$ bisects $\gamma'$ and since $u'_{1,1}\in \partial_{\delta}(G')$, by Lemma \Ref{lab-40}, we get that
$$d'(x'_2,w'_1)\leq \nu_2\delta_{G'}(w'_1)\leq \nu_2 d'(u'_{1,1},w'_1),$$
which, together with (\ref{Add25-1}), shows
\beqq d'(x'_{1,2},u'_{1,1})&<&  2\nu_2 \delta_{G'}(w'_1)+d'(u'_{1,1},w'_1)\\ \nonumber&\leq&(1+2\nu_2)d'(u'_{1,1},w'_1)
\\ \nonumber&\leq&  \frac{2\nu_2+1}{4\nu_2}d'(u'_{1,2},u'_{1,1}), \;\;\;\;\;\;\;(\mbox{by Lemma \ref{Hz-W-Zhou-1}})\eeqq
and thus,
\beqq d'(u'_{1,2},x'_{1,2})\geq d'(u'_{1,2},u'_{1,1})-d'(x'_{1,2},u'_{1,1})>\frac{1}{2}d'(x'_{1,2},u'_{1,1}),\eeqq
which proves the lemma.
\epf

Now, we are ready to get a contradiction of the assumption \eqref{Exter-1}. Since $f$ is $\eta$-quasisymmetric on the $\partial_{\sigma}G$ (by Condition \ref{PropertyB}), we infer from Lemmas \ref{Zhou-W-H-17} and \ref{imp-1}
that $$\frac{1}{2}<\frac{d'(u'_{1,2},x'_{1,2})}{ d'(x'_{1,2},u'_{1,1})}\leq \eta\Big(\frac{\sigma(u_{1,2},x_{1,2})}{\sigma(x_{1,2},u_{1,1})}\Big)\leq\eta\Big(270M_1^{-\frac{1}{6}}\Big).$$
This is impossible because $M_1\geq e^{MM_2^5}$ and $M_2\geq 10\nu_3^4\eta(\nu_3) \big(\eta^{-1}(\frac{1}{5}\nu_2^{-3})\big)^{-1}\geq 10\nu_3^4\big(\eta^{-1}(\frac{1}{2})\big)^{-1}.$
This contradiction implies
that the inequality in (\ref{Exter-1}) is not true either, proving the proposition in this subcase.
\epf

Now, we conclude from \eqref{25-1-13-1} and Proposition \ref{lem-6.1} that the theorem is true.
\epf

%%%%%%%%%%%%%%%%%%%%%%%%%%%%%%%%%%%%%%%%%%%
%%%%%%%%%%%%%%%%%%%%%%%%%%%%%%%%%%%%%%%%%%%
\section{Inner uniformity under Condition \ref{PropertyB}}\label{sec-6}
%%%%%%%%%%%%%%%%%%%%%%%%%%%%%%%%%%%%%%%%%%%
%%%%%%%%%%%%%%%%%%%%%%%%%%%%%%%%%%%%%%%%%%%

The purpose of this section is to complete the proof of Theorem \ref{thm-2.3}. By Proposition \ref{thm-2.1} in Section \ref{sec-3}, still, we have to demonstrate the implication from $(\mathfrak{c})$ to $(\mathfrak{a})$, which is formulated in the following theorem.

\begin{thm}\label{thm-6.1}
 Suppose that $G$ is a bounded domain in $\mathbb{R}^n$ $(n\geq 2)$ and satisfies Condition \ref{PropertyB}.
Then $G$ is $B'$-inner uniform, where $B'$ is a constant depending only on the constant $B_0$ in Theorem \Ref{Thm-HRWZ} and the dimension $n$.
\end{thm}

Before the proof of the theorem, let us recall the concept of a cigar, due to N\"{a}kki and V\"ais\"al\"a \cite{RJ}.

Let $D$ denote a domain in $\IR^n$ $(n\geq 2)$, $u, v\in D$ with $u\not=v$, and let $\alpha\subset D$ be a rectifiable curve connecting $u$ and $v$. For $c\geq 1$, set
$$\mathop{\mathrm{cig}_{\ell}}(\alpha, c)=\bigcup\Big\{\IB\big(w, \frac{1}{c}\rho_{\ell}(w)\big):\; w\in \alpha\backslash \{u, v\}\Big\}$$
and
$$\mathop{\mathrm{cig}_{d}}(\alpha, c)=\bigcup\Big\{\IB\big(w, \frac{1}{c}\rho_{d}(w)\big):\; w\in \alpha\backslash \{u, v\}\Big\},$$
where $$\rho_{\ell}(w)=\min\{\ell(\alpha[u, w]),\; \ell(\alpha[w, v])\}$$ and
$$\rho_{d}(w)=\min\{\diam(\alpha[u, w]),\; \diam(\alpha[w, v])\}.$$

We call $\mathop{\mathrm{cig}_{\ell}}(\alpha, c)$ (resp. $\mathop{\mathrm{cig}_{d}}(\alpha, c)$) the {\it length $c$-cigar} (resp. the {\it diameter $c$-cigar}) with core $\alpha$ joining $u$ and $v$.

A domain $D$ in $\IR^n$ is said to satisfy the {\it condition $c$-$\mathop{\mathrm{cig}_{\ell}}$} (resp. the {\it condition $c$-$\mathop{\mathrm{cig}_{d}}$}) if each pair of distinct points of $D$ can be joined by a $c$-length cigar $\mathop{\mathrm{cig}_{\ell}}(\alpha, c)$ (resp. a $c$-diameter cigar $\mathop{\mathrm{cig}_{d}}(\alpha, c)$). Here, two points in $D$ which are joined by a $\mathop{\mathrm{cig}_{\ell}}(\alpha, c)$ $\mathop{\mathrm{cig}_{d}}(\alpha, c)$) means that $\mathop{\mathrm{cig}_{\ell}}(\alpha, c)\subset D$ (resp. $\mathop{\mathrm{cig}_{d}}(\alpha, c)\subset D$) and $\alpha$ connects these two points.
 Then the following corollary is a direct consequence of the definitions.

 \bcor\label{8-16-1}
Let $u$ and $v$ be two points in $ D$. Then
\ben
\item\label{8-16-1-1}
there is a $c$-diameter cigar $\mathop{\mathrm{cig}_{d}}(\alpha, c)$ connecting $u$ and $v$ if and only if for any $x\in \alpha$, $$\min\{\diam(\alpha[u, x]),\; \diam(\alpha[x, v])\}\leq c\delta_D(x);$$
\item\label{8-16-1-2}
there is a $c$-length cigar $\mathop{\mathrm{cig}_{\ell}}(\alpha, c)$ connecting $u$ and $v$ if and only if for any $x\in \alpha$, $$\min\{\ell(\alpha[u, x]),\; \ell(\alpha[x, v])\}\leq c\delta_D(x).$$
\een
\ecor

The following results are from \cite{RJ}.

\begin{Thm}\label{Vai-1} Suppose that $D$ is a domain in $\mathbb{R}^n$ $(n\geq 2)$.
\ben
\item\label{Vai-4}   $($\cite[Proof of Theorem 2.14]{RJ}$)$
Assume that $D$ satisfies the condition $c$-$\mathop{\mathrm{cig}_{d}}$. Then there is a constant $c_1=c_1(c, n)$ such that for any diameter $c$-cigar $\mathop{\mathrm{cig}_{d}}(\alpha, c)$ connecting two distinct points in $D$, there is a length $c_1$-cigar $\mathop{\mathrm{cig}_{\ell}}(\alpha_1, c_1)$ also connecting these two points which satisfies $$\mathop{\mathrm{cig}_{\ell}}(\alpha_1, c_1)\subset \mathop{\mathrm{cig}_{d}}(\alpha, c).$$

\item$($\cite[Theorem 2.14]{RJ}$)$
$D$ satisfies the condition $c$-$\mathop{\mathrm{cig}_{d}}$ if and only if it satisfies the condition $c'$-$\mathop{\mathrm{cig}_{\ell}}$, where $c'$ depends on $c$ and $n$, and $c$ depends on only $c'$.
\een
\end{Thm}

Now, we are ready to prove the theorem.

\subsection*{Proof of Theorem \ref{thm-6.1}}
Assume that $G$ satisfies Condition \ref{PropertyB}.
Let $x_1$, $x_2$ be two points in $G$. Then there is a curve $\gamma$ in $G$ connecting $x_1$, $x_2$ such that its image $\gamma'=f(\gamma)$ is a quasihyperbolic geodesic in $G'$. Since Theorem \Ref{Thm-HRWZ} and Corollary \ref{8-16-1}\eqref{8-16-1-1} guarantee that $\mathop{\mathrm{cig}_{d}}(\gamma, B_0)\subset G$, it follows from Theorem \Ref{Vai-1}\eqref{Vai-4} that there is a length $B_1$-cigar $\mathop{\mathrm{cig}_{\ell}}(\gamma_1, B_1)$ connecting $x_1$ and $x_2$ such that
\be\label{Vai-5}
\mathop{\mathrm{cig}_{\ell}}(\gamma_1, B_1)\subset \mathop{\mathrm{cig}_{d}}(\gamma, B_0),
\ee where $B_1=B_1(B_0, n)$.
Then Corollary \ref{8-16-1}\eqref{8-16-1-2} ensures that for any $x\in \gamma_1$,
 $$\min\{\ell(\gamma_1[x_1, x]),\; \ell(\gamma_1[x, x_2])\}\leq B_1 \delta_G(x).$$

Let $x_0$ be the point in $\gamma_1$ which bisects $\gamma_1$. Again, \eqref{Vai-5} guarantees that
$$\IB\big(x_0, \frac{1}{2B_1}\ell(\gamma_1)\big)\subset \mathop{\mathrm{cig}_{d}}(\gamma, B_0)\subset \IB(x_1, 2\diam(\gamma)).$$
This leads to
$$\ell(\gamma_1)\leq 4B_1\diam(\gamma).$$
 Now, we conclude from Theorem \ref{sat-11} that
 $$\ell(\gamma_1)\leq 4B_0B_1\sigma(x_1,x_2),$$
where the constant $B_0$ is from Theorem \Ref{Thm-HRWZ}.

So far, we have shown that $\gamma_1$ is $4B_0B_1$-inner uniform in $G$ connecting  $x_1$ and $x_2$. By letting $B'=4B_0B_1$, we know from the arbitrariness of the pair of points $x_1$ and $x_2$ in $G$ that $G$ is $B'$-inner uniform.
\qed

%%%%%%%%%%%%%%%%%%%%%%%%%%%%%%%%%%%%%%%%%%%
%%%%%%%%%%%%%%%%%%%%%%%%%%%%%%%%%%%%%%%%%%%
\section{Proofs of the implications $(\mathfrak{ii})$ $\Longrightarrow$ $(\mathfrak{iii})$  $\Longrightarrow$ $(\mathfrak{i})$ in Theorem \ref{thm-1.1}}\label{sec-7}
%%%%%%%%%%%%%%%%%%%%%%%%%%%%%%%%%%%%%%%%%%%
%%%%%%%%%%%%%%%%%%%%%%%%%%%%%%%%%%%%%%%%%%%

In this section, we complete the proof of Theorem \ref{thm-1.1} by showing the implications $(\mathfrak{ii})$ $\Longrightarrow$ $(\mathfrak{iii})$ $\Longrightarrow$ $(\mathfrak{i})$.

\subsection{Proof of the implication $(\mathfrak{ii})$ $\Longrightarrow$ $(\mathfrak{iii})$}

Before the proof, let us introduce a concept.

Let $c\geq 1$. A  non-complete metric space $(X,d)$ is {\it $c$-locally externally connected}, or $c$-{\it LEC}, with respect to the metric $d$ if for all points $x\in X$ and all $r\in (0,\delta_X(x)/c)$,
every pair of points in $X\setminus \overline{\mathbb{B}}_d(x,r)$ can be joined by a rectifiable curve in  $X\setminus \overline{\mathbb{B}}_d(x,r/c)$.

It is not difficult to see that {\it LLC$_2$} implies {\it LEC}, and every Euclidean domain in $\IR^n$ with $n\geq 2$ is $c$-{\it LEC} for each $c>1$ with respect to the Euclidean metric, and its corresponding inner metric as well.
The following useful theorem concerning the {\it LEC} property is due to Buckley and Herron.

\begin{Thm}\label{bhthm4.2} $($\cite[Theorem 4.2]{BH}$)$ Suppose that $(X, d)$ is rectifiably connected, locally compact and non-complete. If it is uniform and $LEC$, then it must be {\it LLC}.
\end{Thm}

Now, we are ready to prove the implication $(\mathfrak{ii})$ $\Longrightarrow$ $(\mathfrak{iii})$. Assume that $G$ satisfies the condition $(\mathfrak{ii})$. Then by the proved equivalence of the conditions $(\mathfrak{i})$ and $(\mathfrak{ii})$ in Theorem \ref{thm-1.1}, or Theorem \ref{thm-2.3}, we know that $G$ is inner uniform, and therefore, $(G,\sigma)$ is uniform. Since $(G,\sigma)$ is $LEC$, we know from Theorem \Ref{bhthm4.2} that $(G,\sigma)$ is {\it LLC}. This shows that $G$ is {\it LLC} with respect to the inner metric, and hence, the implication is proved.\qed

\subsection{Proof of the implication $(\mathfrak{iii})$ $\Longrightarrow$ $(\mathfrak{i})$}

For the proof, assume that $G$ is Gromov $\delta$-hyperbolic and $c$-{\it LLC} with respect to the inner metric $\sigma$, where $\delta\geq 0$ and $c\geq 1$, i.e., $G$ satisfies the condition $(\mathfrak{iii})$ in Theorem \ref{thm-1.1}.
Then Theorem \Ref{BB03} shows that
$G$ satisfies the $c_{gh}$-Gehring-Hayman condition and the $c_{bs}$-ball-separation condition with $c_{gh}\geq 1$ and $c_{bs}\geq 1$.

To prove the inner uniformity of $G$, let $x$ and $y$ be two distinct points in $G$, and let $\gamma=[x, y]_k$, a quasihyperbolic geodesic in $G$ connecting $x$ and $y$. Since $G$ satisfies the $c_{gh}$-Gehring-Hayman condition, we know that
\be\label{zhou-1}
\ell(\gamma)\leq c_{gh}\sigma(x,y).
\ee

In the following, we show that for all $z\in \gamma$,
\be\label{zhou-2}
\min\{\ell(\gamma[x,z]),\ell(\gamma[z,y])\}\leq cc_{gh}c_{bs} \delta_G(z).
\ee

For this, let $z\in \gamma$ and
$$t(z)=\frac{\min\{\ell(\gamma[x,z]),\ell(\gamma[z,y])\}}{\delta_G(z)}.$$

On the one hand, by the $c_{gh}$-Gehring-Hayman condition, we see that
$$\ell(\gamma[x,z])\leq c_{gh}\sigma(x,z)\;\;{\rm and}\;\;\ell(\gamma[y,z])\leq c_{gh}\sigma(y,z).$$
These ensure that
$$\min\{\sigma(x,z),\sigma(y,z)\}\geq \frac{1}{c_{gh}}t(z)\delta_G(z),$$ from which we know that
$$x,\,y\notin \mathbb{B}_\sigma\Big(z,\frac{1}{c_{gh}}t(z)\delta_G(z)\Big).$$
Because $G$ is $c$-{\it LLC} with respect to the inner metric $\sigma$, it follows that there is a curve $\alpha$ in $G$ connecting $x$ and $y$ such that
$$\alpha\subset G\setminus \overline{\mathbb{B}}_\sigma\Big(z,\frac{1}{cc_{gh}}t(z)\delta_G(z)\Big).$$
Thus we have
\be\label{z-1}\dist_\sigma(z, \alpha)\geq \frac{1}{cc_{gh}}t(z)\delta_G(z).\ee

On the other hand, since $\gamma$ is a quasihyperbolic geodesic and because $G$ satisfies the $c_{bs}$-ball-separation condition, we know that
$$\alpha\cap \mathbb{B}_\sigma(z,c_{bs}\delta_G(z))\not=\emptyset.$$
By combining with (\ref{z-1}), we obtain
$$\frac{1}{cc_{gh}}t(z)\delta_G(z)\leq \dist_\sigma(z, \alpha) \leq c_{bs}\delta_G(z),$$
which implies
$$t(z)\leq cc_{gh}c_{bs}.$$ Because this holds for all $z$ in $\gamma$, we see that the inequality \eqref{zhou-2} is true.
Now, we conclude from \eqref{zhou-1} that $G$ is $cc_{gh}c_{bs}$-inner uniform. \qed

\subsection*{Acknowledgments}
We are grateful to the anonymous referee for the constructive suggestions, which have helped us to significantly improve the quality of this article. 

The first author (Manzi Huang) was partly supported by NSF of China under the number 12371071, and the second author (Antti Rasila) was partly supported by NSF of China under the number 11971124, NSF of Guangdong Province under the numbers 2021A1515010326 and 2024A1515010467, and Li Ka Shing Foundation STU-GTIIT Joint-Research Grant under the number 2024LKSFG06; the third author (Xiantao Wang) was partly supported by NSF of China under the number  12071121, and the fourth author (Qingshan Zhou) was partly supported by NSF of China under the number 11901090.

\end{document}